\documentclass[fleqn,a4paper,10pt]{amsart}

\usepackage{graphicx}
\usepackage{tikz}
\usetikzlibrary{shapes,decorations.pathmorphing}
\usepackage{amsmath}
\usepackage{amssymb,amsthm}
\usepackage{comment}
\usepackage{hyperref}
\usepackage{enumerate}

\theoremstyle{definition}

\title[A universal law for Vorono\"\i\ cell volumes in infinitely large maps]{A universal law for Vorono\"\i\ cell volumes in infinitely large maps}
\author{Emmanuel Guitter}
\address{Institut de physique th\'eorique, Universit\'e Paris Saclay, CEA, CNRS, F-91191 Gif-sur-Yvette}
\email{emmanuel.guitter@ipht.fr}

\begin{document}
\maketitle

\begin{abstract}
We discuss the volume of Vorono\"\i\ cells defined by two marked vertices picked randomly at a fixed given mutual distance $2s$ in random planar quadrangulations. 
We consider the regime where the mutual distance $2s$ is kept finite while the total volume of the quadrangulation tends to infinity. In
this regime, exactly one of the Vorono\"\i\ cells keeps a finite volume, which scales as $s^4$ for large $s$. 
We analyze the universal probability distribution 
of this, properly rescaled, finite volume and present an explicit formula for its Laplace transform.  
 \end{abstract}

\section{Introduction}
\label{sec:introduction}

In a recent paper \cite{G17a}, we analyzed the volume distribution of Vorono\"\i\ cells for some families of random bi-pointed planar maps.
Recall that a planar map is a connected graph embedded in the sphere: it is bi-pointed if it has two marked distinct vertices.
These marked vertices allow us to partition the map into two Vorono\"\i\ cells, where each cell corresponds, so to say, to the part
of the map closer to one marked vertex than to the other. The \emph{volume} of, say the second Vorono\"\i\ cell (that centered
around the second marked vertex) is then a finite fraction $\phi$ of the total volume of the map, with $0\leq \phi\leq 1$, while the first cell 
clearly spans the complementary fraction $1-\phi$. The main result proven in \cite{G17a} is that, for several families of random bi-pointed maps
with a \emph{fixed total volume}, and in the limit where this volume becomes infinitely large, \emph{the law for the fraction $\phi$ of the total volume
spanned by the second Vorono\"\i\ cell is uniform in the interval $[0,1]$}, a property conjectured by Chapuy in \cite{GC16}
among other more general conjectures.
Here it is important to stress that the above result holds when the two marked vertices are chosen \emph{uniformly at random} in the map.
In particular, 
their mutual distance is left arbitrary\footnote{For one family on maps considered in \cite{G17a}, it was assumed for convenience that 
the mutual distance be even, but lifting this constraint has no influence on the obtained result.}.

This paper deals on the contrary with Vorono\"\i\ cells within random bi-pointed maps where the two marked vertices are picked randomly 
\emph{at a fixed given mutual distance}. Considering again the limit of maps with an infinitely large volume and \emph{keeping the 
(fixed) mutual distance between the marked vertices finite}, we find that only one of the Vorono\"\i\ cells
becomes infinitely large while the volume of the other remains finite. In particular, the fraction of the total volume spanned by 
this latter cell tends to $0$ while that of the infinite cell tends to $1$. In other words, having imposed a fixed finite mutual distance between the marked vertices drastically modifies the law for 
the fraction $\phi$ which is now concentrated at $\phi=0$ if it is precisely the second Vorono\"\i\ cell which remains finite or 
at $\phi=1$ if this second cell becomes infinite.  

In this regime of fixed mutual distance, a good measure of the  Vorono\"\i\ cell extent is now provided by  the \emph{volume of that of the two Vorono\"\i\ cells which remains finite}. The main goal of this paper is to compute the law for this finite volume, in particular in a universal
regime where the mutual distance, although kept finite, is large.

The paper is organized as follows: we first introduce in Section 2 the family of bi-pointed maps that will shall study (i.e.\ bi-pointed quadrangulations),
define the volumes of the associated Vorono\"\i\ cells and introduce some generating function with some control on these volumes
(Section 2.1). We then discuss the \emph{scaling function} which captures the properties of this generating function
in some particular scaling regime (Section 2.2), and whose knowledge is the key of the subsequent calculations.  Section 3 is devoted to our
analysis of 
Vorono\"\i\ cell volumes in the regime of interest in this paper, namely when the maps become infinitely large and the mutual distance between the marked vertices remains finite.
We first analyze (Section 3.1) the law for the fraction $\phi$ of the total volume of the map spanned by the second Vorono\"\i\ cell and
show, as announced above, that it is evenly concentrated at  $\phi=0$ or $\phi=1$. We then analyze (Section 3.2) map configurations for which the
volume of the second Vorono\"\i\ cell remains finite and show how to obtain, from 
the simple knowledge of the scaling function introduced above, the law for this (properly rescaled) volume when the mutual distance becomes large. This leads to an explicit universal expression (Section 3.3) for the 
probability distribution of the finite Vorono\"\i\ cell volume (in practice for its Laplace transform), whose properties are discussed in details. Section 4 proposes 
an instructive comparison of our result with that, much simpler, obtained for Vorono\"\i\ cells within bi-pointed random trees. Section
5 discusses the case of asymmetric Vorono\"\i\ cells where some explicit bias in the evaluation of distances is introduced.
Our conclusions are gathered in Section 6. A few technical details, as well as explicit but heavy intermediate expressions, are given
in various appendices. 

\section{Vorono\"\i\ cells in bi-pointed maps}
\label{sec:voronoi}
\subsection{A generating function for bi-pointed maps with a control on their Vorono\"\i\ cell volumes}

The objects under study in this paper are \emph{bi-pointed planar quadrangulations}, namely planar maps whose all faces have degree $4$, and with two marked distinct vertices.
We moreover demand that these vertices, distinguished as $v_1$ and $v_2$, be at some even graph distance $d(v_1,v_2)$, namely
\begin{equation*}
d(v_1,v_2)=2s
\end{equation*}
for some fixed given integer $s\geq 1$. 
Given $v_1$ and $v_2$, the corresponding two Vorono\"\i\ cells are obtained via some splitting of the map
into two domains which, so to say, regroup vertices which are closer to one marked vertex than to the other. 
As discussed in details in \cite{G17a}, a canonical way to perform this splitting consists in applying the well-know Miermont bijection \cite{Miermont2009}
which 
transforms a bi-pointed planar quadrangulation into a so-called planar \emph{iso-labelled two-face map (i-l.2.f.m)}, namely a planar map with exactly two faces, distinguished as 
$f_1$ and $f_2$ and with vertices labelled by positive integers satisfying:
\begin{enumerate}[$(\hbox{L}_1)$]
\item{labels on adjacent vertices differ by $0$ or $\pm1$;}
\item{the minimum label for the set of vertices incident to $f_1$ is $1$;}
\item{the minimum label for the set of vertices incident to $f_2$ is $1$.}
\end{enumerate}
As recalled in \cite{G17a}, the  Miermont bijection provides a one-to-one correspondence between bi-pointed planar quadrangulations and planar i-l.2.f.m, the labels of the vertices corresponding precisely to their distance to the closest marked vertex in the quadrangulation. More interestingly, by drawing the original quadrangulation on top of its image, the two faces $f_1$ and $f_2$ define {\it de facto} two domains in the quadrangulation which are perfect realizations of the desired two Vorono\"\i\ cells as, by construction, each of these domains regroups vertices closer to one marked vertex. 
Since faces of the quadrangulation are, under the Miermont bijection, in correspondence with edges of the i-l.2.f.m, 
the volume ($=$ number of faces) of a given cell
in the quadrangulation is measured by \emph{half the number of edge sides} incident to the corresponding face in the i-l.2.f.m. 
Note that this volume is in 
general some half-integer since a number of faces of the quadrangulation may be shared by the two cells (see \cite{G17a} for details). 
To be precise, an i-l.2.f.m is made of a simple closed loop $\mathcal{L}$ separating its two faces $f_1$ and $f_2$\footnote{This loop
is simply formed by the cyclic sequence of edges incident to both faces.} together with a number of subtrees attached to vertices along $\mathcal{L}$, possibly on each side of the loop. If we call $e_1$ and $e_2$ the total number of edges for subtrees in the face $f_1$ and $f_2$ respectively, and $e$ the length ($=$ number of edges) of the loop $\mathcal{L}$, the volumes $n_1$ and $n_2$ of the Vorono\"\i\ cells are respectively
\begin{equation*}
n_i=e_i+\frac{e}{2}, \qquad i=1,2,
\end{equation*}
for a total volume
\begin{equation*}
N= n_1+n_2= e_1+e_2+e\ .
\end{equation*}
Finally, the requirement that $d(v_1,v_2)=2s$ translates into the following fourth label constraint:
\begin{enumerate}[$(\hbox{L}_4)$]
\item{the minimum label for the set of vertices incident to $\mathcal{L}$ is $s$.}
\end{enumerate}

Having defined Vorono\"\i\ cells, we may control their volume by considering the \emph{generating function} $F(s,g,h)$ of bi-pointed planar quadrangulations where $d(v_1,v_2)=2s$, with a weight 
\begin{equation*}
g^{n_1}\, h^{n_2}\ .
\end{equation*}
From the Miermont bijection and the associated canonical construction of Vorono\"\i\ cells, $F(s,g,h)$ is also the generating function of i-l.2.f.m satisfying the extra requirement $(\hbox{L}_4)$ with a weight 
\begin{equation*}
g^{e_1}\, h^{e_2}\ (\sqrt{g\, h})^e\ .
\end{equation*}
As such, $F(s,g,h)$ may, via some appropriate decomposition of the i.l.2.f.m, be written as (see \cite{G17a})
\begin{equation}
F(s,g,h)= \Delta_s\Delta_t \log(X_{s,t}(g,h))\Big\vert_{t=s} =\log\left(\frac{X_{s,s}(g,h)X_{s-1,s-1}(g,h)}{X_{s-1,s}(g,h)X_{s,s-1}(g,h)}\right)
\label{eq:FsX}
\end{equation}
(here $\Delta_s$ is the finite difference operator $\Delta_s f(s) \equiv f(s)-f(s-1)$), 
where $X_{s,t}(g,h)$ is some generating function for appropriate \emph{chains} of labelled trees (which correspond to appropriate open
sequences of edges with subtrees attached on either side of the incident vertices). Without entering into details, it is enough for the scope of
this paper to know that the generating function $X_{s,t}(g,h)$  is entirely determined\footnote{This relation fully determines $X_{s,t}(g,h)$ for all $s,t\geq 0$ order by order in $g$ and $h$, i.e.\ 
$X_{s,t}(\rho g,\rho h)$ is fully determined order by order in $\rho$.} by the relation (obtained by a simple splitting of the chains)
\begin{equation}
X_{s,t}(g,h)=1+\sqrt{g\, h}\, R_s(g)R_t(h)X_{s,t}(g,h)\left(1+\sqrt{g\, h}\, R_{s+1}(g)R_{t+1}(h)X_{s+1,t+1}(g,h)\right)
\label{eq:eqforXst}
\end{equation}
for $s,t\geq 0$, where 
the quantity $R_s(g)$ (as well as its analog $R_t(h)$) is a well known generating function for appropriate 
labelled trees. It is given explicitly by 
\begin{equation}
R_s(g)=\frac{1+4x+x^2}{1+x+x^2}\frac{(1-x^s)(1-x^{s+3})}{(1-x^{s+1})(1-x^{s+2})}
\quad \hbox{for}\ g= x\frac{1+x+x^2}{(1+4x+x^2)^2}\ ,
\label{eq:eqforRs}
\end{equation}
where $x$ is taken in the range $0\leq x\leq 1$ and parametrizes $g$ (in the range $0\leq g\leq 1/12$ for a proper convergence of the generating function).
For $h=g$, the solution of \eqref{eq:eqforXst} can be made explicit and reads
\begin{equation}
X_{s,t}(g,g)=\frac{(1-x^3)(1-x^{s+1})(1-x^{t+1})(1-x^{s+t+3})}{(1-x)(1-x^{s+3})(1-x^{t+3})(1-x^{s+t+1})}\ .
\label{eq:exactXgg}
\end{equation}
Unfortunately, no such explicit expression is known for $X_{s,t}(g,h)$ when $h\neq g$ and the relation \eqref{eq:FsX} might thus appear of no practical use at a first glance.
As discussed in \cite{G17a}, this is not quite true as we may recourse to appropriate scaling limits of all the above generating functions to extract explicit statistics on Vorono\"\i\ cell volumes
in a limit where the maps become (infinitely) large. Let us now discuss this point.
 
\subsection{The associated scaling function}
The limit of large quadrangulations (i.e.\ with a large number $N$ of faces) is captured by the singularity of $F(s,g,h)$ whenever $g$ or $h$ tends
toward its critical value $1/12$. As we shall see, in all cases of interest, this singularity may be analyzed by setting
\begin{equation}
g=G(a,\epsilon)\ ,\qquad h=G(b,\epsilon)\ , \quad \hbox{where}\ \  G(c,\epsilon)\equiv \frac{1}{12}\left(1-\frac{c^4}{36}\epsilon^4\right)\ ,
\label{eq:ghscal}
\end{equation}
and letting $\epsilon$ tend to $0$. In this limit, we have for instance the following expansion for the quantity $x$ parametrizing $g$
in \eqref{eq:eqforRs}:
\begin{equation*}
x=1-a\, \epsilon +\frac{a^2 \epsilon ^2}{2}-\frac{5\, a^3 \epsilon ^3}{24}+\frac{a^4 \epsilon ^4}{12}-\frac{13\, a^5 \epsilon ^5}{384}+\frac{a^6 \epsilon ^6}{72}-\frac{157\, a^7 \epsilon ^7}{27648}+\frac{a^8 \epsilon
   ^8}{432}+O(\epsilon ^{9})\ ,
\end{equation*}
so that, for $h=g$ (i.e.\ $b=a$), we easily get from the exact expression \eqref{eq:exactXgg} of $X_{s,t}(g,g)$ the expansion
\begin{equation}
\begin{split}
F(s,g,g) &=\log \left(\frac{s^2 (2 s+3)}{(s+1)^2 (2 s-1)}\right)-\frac{(2 s+1)a^4 \, \epsilon ^4}{60} 
\\
&\hskip 3.cm +\frac{(2 s+1) \left(10 s^2+10 s+1\right) a^6 \, \epsilon ^6}{1890}+O(\epsilon ^8)\ .\\
   \end{split}
   \label{eq:expDelta}
\end{equation}
Since $a^4\epsilon^4=36 (1-12g)$ is regular when $g\to 1/12$, the most singular part of this generating function is given by 
 \begin{equation*}
 \frac{(2 s+1) \left(10 s^2+10 s+1\right) a^6 \, \epsilon ^6}{1890}=  \frac{4\,  (2 s+1) \left(10 s^2+10 s+1\right)}{35}\, (1-12g)^{3/2}
 \end{equation*}
 and we thus deduce that the \emph{number} $ F_N(s)$ of bi-pointed planar quadrangulations with $N$ faces and with their two marked vertices at distance $2s$ behaves at large $N$ as 
 \begin{equation}
\hskip -10pt F_N(s)\equiv  [g^N]F(s,g,g)\underset{N \to \infty}{\sim} \frac{3}{4}  \frac{12^{N}}{\sqrt{\pi} N^{5/2}} \mathfrak{f}_3(s)\ ,\qquad 
\mathfrak{f}_3(s)=\frac{4\,  (2 s+1) \left(10 s^2+10 s+1\right)}{35}\ .
\label{eq:number}
 \end{equation}
 When $s$ itself becomes large, this number behaves as
 \begin{equation}
\frac{3}{4}\, \frac{12^{N}}{\sqrt{\pi} N^{5/2}}  \times \frac{16}{7}s^3\ .
   \label{eq:larges}
   \end{equation}
 Note that this later estimate assumes that $N$ becomes \emph{first} arbitrarily large with a value of $s$ remaining finite, and \emph{only then} is $s$ set to be large. This order of limits
 corresponds to what is usually called the \emph{local limit}. In particular, $N$ and $s$ do not scale with each other.
 \vskip .3cm  
Now it is interesting to note that getting this last result \eqref{eq:larges} does not require the full knowledge of $F(s,g,g)$ and may be obtained upon using
instead  some simpler \emph{scaling function} which captures the behavior of $F(s,g,g)$ 
in a particular \emph{scaling regime}. Consider indeed the generating function $X_{s,t}(g,g)$ in a regime where $g\to 1/12$ as above by letting $\epsilon \to 0$ in \eqref{eq:ghscal}, but where we let \emph{simultaneously} $s$ and $t$ become large
upon setting
\begin{equation*}
s=\frac{S}{\epsilon}\ , \qquad t=\frac{T}{\epsilon}\ ,
\end{equation*}
with $S$ and $T$ kept finite. In this scaling regime, we have the expansion
\begin{equation*}
X_{\left\lfloor S/\epsilon\right\rfloor,\left\lfloor T/\epsilon\right\rfloor}(g,g)  =3+x(S,T,a)\ \epsilon+O(\epsilon^2)\ ,
\end{equation*}
where the function $x(S,T,a)$ is given explicitly from \eqref{eq:exactXgg} by
\begin{equation*}
x(S,T,a) =-3\, a-\frac{6 a \left(e^{-a S}+e^{-a T}-3 e^{-a (S+T)}+e^{-2a (S+T)}\right)}{\left(1-e^{-a S}\right) \left(1-e^{-a T}\right) \left(1-e^{-a (S+T)}\right)}\ .
\end{equation*}  
This in turn implies the expansion
\begin{equation*}
\begin{split}
\Delta_s\Delta_t \log\left(X_{\left\lfloor S/\epsilon\right\rfloor,\left\lfloor T/\epsilon\right\rfloor}(g,g) \right) \Big\vert_{T=S}&=
 \partial_S\partial_T \log\left(3+x (S,T,a)\, \epsilon \right)\Big\vert_{T=S}\times \epsilon^2+O(\epsilon^4)\\ &=
 \frac{1}{3}\ \partial_S\partial_T x(S,T,a)\Big\vert_{T=S}\times \epsilon^3 +O(\epsilon^4)\\
\end{split}
\end{equation*}
which yields
\begin{equation*}
F\left(\left\lfloor S/\epsilon\right\rfloor,g,g \right)=
\mathcal{F}(S,a)\, \epsilon^3+O(\epsilon^4)\ ,
\end{equation*}
where the \emph{scaling function} $\mathcal{F}(S,a)$ associated with $F[s,g,g]$ reads explicitly
\begin{equation}
\begin{split}
\mathcal{F}(S,a)&=\frac{1}{3}\ \partial_S\partial_T x(S,T,a)\Big\vert_{T=S} \\
& = \frac{2 \, a^3\, e^{-2 a S} \left(1+e^{-2 a S}\right)}{\left(1-e^{-2 a S}\right)^3} \\
&= \left(\frac{1}{2\, S^3}-\frac{a^4 S}{30}+\frac{2 a^6 S^3}{189}+O(S^5)\right)\ .\\
\end{split}
\label{eq:scalF}
\end{equation}
A crucial remark is that we recognize in this latter \emph{small $S$ expansion} of $\mathcal{F}(S,a)$ the \emph{large $s$ leading behavior}\footnote{In particular, we have the large $s$ expansion:
$\log \left(\frac{s^2 (2 s+3)}{(s+1)^2 (2 s-1)}\right)= \frac{1}{2\, s^3}+O\!\left(\frac{1}{s^4}\right)$\ .}
of the coefficients in the expansion \eqref{eq:expDelta} for $F[s,g,g]$ in the local limit. For instance, the large $s$ 
behavior of the singular term (proportional to $\epsilon^6$) in \eqref{eq:expDelta} is given by 
\begin{equation*}
\frac{(2 s+1) \left(10 s^2+10 s+1\right) a^6}{1890}\underset{s \to \infty}{\sim} \frac{2}{189} s^3\, a^6 \, = s^3\, \times [S^3] \mathcal{F}(S,a)\ .
\end{equation*}
For $a=\sqrt{6}$ (in which case we have the direct identification $\epsilon^6=(1-12g)^{3/2}$),  the left hand side is precisely the coefficient $\mathfrak{f}_3(s)$ in \eqref{eq:number}, so that the result \eqref{eq:larges} may thus be read off \emph{directly on the expression of the scaling function $\mathcal{F}(S,a)$} via
\begin{equation}
\mathfrak{f}_3(s)\underset{s \to \infty}{\sim}  s^3\, \times [S^3] \mathcal{F}(S,\sqrt{6})=  s^3\, \times \frac{2}{189} (\sqrt{6})^6 = \frac{16}{7}\, s^3\ ,
\label{eq:largesbis}
\end{equation}
without recourse to the explicit knowledge of
the full generating function $F[s,g,g]$.

The origin of this \emph{``scaling correspondence"}, which connects the local limit at large $s$ to the scaling limit at small $S$ is explained in details
in the next section. This correspondence is in fact a general property and can be applied in the situation where $h\neq g$. It therefore allows us to access the large $s$ limit of the
large $N$ asymptotics of $[g^N]F(s,g,h)$ (again sending $N\to \infty$ first) from the simple knowledge of the scaling function associated 
with $F(s,g,h)$.  
\vskip .3cm
As of now, let
us already fix our notations for scaling functions when $g$ and $h$ are arbitrary: parametrizing $g$ and $h$ as in \eqref{eq:ghscal} above, we have when $\epsilon\to 0$ the expansion
\begin{equation*}
X_{\left\lfloor S/\epsilon\right\rfloor,\left\lfloor T/\epsilon\right\rfloor}(g,h)  =3+x(S,T,a,b)\ \epsilon+O(\epsilon^2)
\end{equation*}
with a scaling function $x(S,T,a,b)$ which, from \eqref{eq:eqforXst} expanded at lowest non-trivial order in $\epsilon$, is solution of the non-linear partial
differential equation 
\begin{equation}
2 \big(x(S,T,a,b)\big)^2+6 \big(\partial_Sx(S,T,a,b)+\partial_Tx(S,T,a,b)\big)+27 \big(r(S,a)+r(T,b)\big)=0\ .
\label{eq:eqforxstab}
\end{equation}
Here $r(S,a)$ is the first non-trivial term in the small $\epsilon$ expansion of $R_{\left\lfloor S/\epsilon\right\rfloor}(g)$, namely, from its explicit expression \eqref{eq:eqforRs},
\begin{equation}
R_{\left\lfloor S/\epsilon\right\rfloor}(g) =2+r(S,a)\ \epsilon^2+O(\epsilon^3)\ , \qquad \ r(S,a)=-\frac{a^2 \left(1+10 e^{-a S}+e^{-2 a S}\right)}{3 \left(1-e^{-a S}\right)^2}\ .
 \end{equation}
As for $F(s,g,h)$, we may now use \eqref{eq:FsX} to relate the associated scaling function $\mathcal{F}(S,a,b)$ to $x(S,T,a,b)$, namely
\begin{equation}
\begin{split}
F(\left\lfloor S/\epsilon\right\rfloor,g,h) &  = \mathcal{F}(S,a,b)\, \epsilon^3+O(\epsilon^4)\\
\hbox{where} & \quad \mathcal{F}(S,a,b) = \frac{1}{3}\ \partial_S\partial_T x(S,T,a,b)\Big\vert_{T=S}\ .\\
\end{split}
\label{eq:xtoF}
\end{equation}
Scaling functions are in general much simpler than the associated full generating functions. In particular, although we have no formula for $X_{s,t}(g,h)$ for arbitrary $g$ and $h$, an explicit expression for $x(S,T,a,b)$ is known for arbitrary $a$ and $b$, 
as first obtained in \cite{G17a} upon solving \eqref{eq:eqforxstab} with appropriate boundary conditions. We may thus recourse to this result to get an explicit expression for the scaling function $\mathcal{F}(S,a,b)$
itself via \eqref{eq:xtoF}. The corresponding formula is quite heavy and its form is not quite illuminating. Still, we display it
in Appendix A for completeness (the reader may refer to this expression to check the various limits and expansions of $\mathcal{F}(S,a,b)$ 
displayed hereafter in the paper).

As opposed to $F(s,g,h)$, the scaling function $\mathcal{F}(S,a,b)$ is thus \emph{known exactly} and we will now show in details how to
use the scaling correspondence to deduce from its small $S$ expansion the large $s$ limit of the
large $N$ asymptotics of $[g^N]F(s,g,h)$ and control the volume of, say, the second Vorono\"\i\ cell in large quadrangulations, by some appropriate choice of $h$.

\section{Infinitely large maps with two vertices at finite distance}

This section is devoted to estimating the law for the volumes spanned by the Vorono\"\i\ cells in bi-pointed quadrangulations whose total volume $N$
($=$ number of faces) tends to infinity. Calling $n_1$ and $n_2$ the two Vorono\"\i\ cell volumes, we have $n_1+n_2=N$ so it is enough to control
one of two volumes, say $n_2$. Here the distance $2s$ between the two marked vertices is kept finite (possibly large) when $N\to \infty$.

\subsection{The law for the proportion of the total volume spanned by one Vorono\"\i\ cell}
For large $N$ and finite $s$, the first natural way to measure $n_2$ is to express it in units of $N$, i.e.\ consider the \emph{proportion}
\begin{equation*}
\phi\equiv \frac{n_2}{N}
\end{equation*}
of the total volume spanned by the second Vorono\"\i\ cell. We have of course $0\leq \phi\leq 1$ and the large $N$ asymptotic probability law $\mathcal{P}_s(\phi)$ for
$\phi$ may be obtained from $F(s,g,h)$ via
\begin{equation*}
\int_{0}^{1}d\phi\, \mathcal{P}_s(\phi)\, e^{\mu \, \phi}= \lim_{N\to \infty} \frac{[g^N]F(s,g,g\, e^{\mu/N})}{[g^N]F(s,g,g)}
\end{equation*}
since $g^{n_1}(g\, e^{\mu/N})^{n_2}=g^N\, e^{\mu \, \phi}$. 

From the scaling correspondence, the large $N$ asymptotics 
of $[g^N]F(s,g,g\, e^{\mu/N})$ is, at large $s$, encoded in the small $S$ expansion of the scaling function $\mathcal{F}(S,a,b)$ for some appropriate $b\equiv b(a,\mu)$.
The right hand side of the above equality may thus be computed explicitly at large $s$ from the knowledge of $\mathcal{F}(S,a,b)$. 
This computation, together with the precise correspondence between the large $s$ local limit and the small $S$ scaling limit, is discussed in details in Appendix B.
We decided however not to develop the calculation here since the resulting law is in fact trivial. As might have been guessed by the reader,
we indeed find
\begin{equation}
\int_{0}^{1}d\phi\, \mathcal{P}_s(\phi)\, e^{\mu \, \phi} \underset{s \to \infty}{\sim}  \frac{1}{2} (1+e^\mu)
\label{eq:trivial}
\end{equation}
or equivalently
\begin{equation*}
 \mathcal{P}_s(\phi) \underset{s \to \infty}{\sim}  \frac{1}{2} (\delta(\phi)+\delta(\phi-1))\ .
\end{equation*}
This result simply states that, for $N\to \infty$  and $s$ finite large, only one of the Vorono\"\i\ cells has a volume of order $N$  with, by symmetry,
\begin{equation}
\left\{
\begin{split}
& n_1=N-o(N)\ , \quad n_2=o(N)\ \quad \hbox{with probability}\ \frac{1}{2}\ ,\\
& n_1=o(N)\ , \quad n_2=N-o(N)\ \quad \hbox{with probability}\ \frac{1}{2}\ .\\
\end{split}
\right.
\label{eq:probahalf}
\end{equation}
The main purpose of this paper, discussed in the following sections,  is precisely to characterize the volume of the Vorono\"\i\ cell which is an $o(N)$. 
As we shall see, the volume of this Vorono\"\i\ cell remains actually finite and scales as $s^4$ when $s$ becomes large.

\subsection{Infinitely large maps with a finite Vorono\"\i\ cell}
This section and the next one present our main result, namely the law for the (properly rescaled) volume of the Vorono\"\i\ cell which is not of order $N$ when $N\to\infty$. More precisely, we will concentrate here on map configurations for which the total volume $N=n_1+n_2$ tends to infinity but 
\emph{the volume $n_2$ is kept finite}. We will then verify {\it a posteriori} that the number of these configurations represents $1/2$ of the total
number of bi-pointed maps whenever $s$ is large. This will {\it de facto} prove that the configurations for which $n_2=o(N)$ in \eqref{eq:probahalf}
are in fact, with probability $1$, configurations for which $n_2$ is finite.

Let us denote by
\begin{equation*}
F_{n_1,n_2}(s)\equiv [g^{n_1}h^{n_2}]F(s,g,h)
\end{equation*}
the number of planar bi-pointed quadrangulations with fixed given values of $n_1$ and $n_2$.  In the limit $N=n_1+n_2\to \infty$ with a 
fixed finite $n_2$, this number may be estimated from the leading singularity of $F(s,g,h)$ when 
$g\to 1/12$ for a fixed value of $h< 1/12$
(see \cite{G16c} for a detailed argument of a fully similar estimate in the context of hull volumes). We have indeed
\begin{equation*}
F_{N-n_2,n_2}(s)\underset{N \to \infty}{\sim} \frac{3}{4}  \frac{12^{N}}{\sqrt{\pi} N^{5/2}}\times 12^{-n_2}\, [h^{n_2}] \mathfrak{f}_3(s,h)\ ,
\end{equation*}
where $\mathfrak{f}_3(s,h)$ is the coefficient of the leading singularity of $F(s,g,h)$ when $g\to 1/12$ \emph{at fixed $h$}, 
hence is obtained from the expansion\footnote{\label{note1}The precise form of this expansion is dictated by the similar explicit expansion \eqref{eq:expDelta} for $F(s,g,g)$. 
In particular, the absence of singular term $\propto (1-12g)^{1/2}$ is imposed by the fact that such a term, if present, would imply that 
$F_{N-n_2,n_2}(s)$  be of order $\hbox{const.}\times 12^N/ N^{3/2}$  at large $N$ while this quantity is clearly bounded
by $[g^N]F[s,g,g]$ which, as we have seen, is of order $\hbox{const'.}\times   12^N/ N^{5/2}$ only.
}
\begin{equation*}
F(s,g,h)=\mathfrak{f}_0(s,h)+\mathfrak{f}_2(s,h)\, (1-12\, g)+\mathfrak{f}_3(s,h)\, (1-12\, g)^{3/2}+O\left((1-12\, g)^2\right)\ .
\end{equation*}
Upon normalizing by the total number of bi-pointed maps $F_N(s)$ with fixed $N$ and $s$, whose asymptotic behavior is given by \eqref{eq:number}, 
we deduce the $N\to\infty$ limiting probability $\mathfrak{p}_s(n_2)$ that the second Vorono\"\i\ cell has volume $n_2$:
\begin{equation*}
\mathfrak{p}_s(n_2)=\lim_{N\to \infty} \frac{F_{N-n_2,n_2}(s)}{F_N(s)} = \frac{1}{\mathfrak{f}_3(s)}\, 12^{-n_2}\, [h^{n_2}] \mathfrak{f}_3(s,h)\ .
\end{equation*}
This probability for arbitrary finite $n_2$ may be encoded in the generating function
\begin{equation}
\sum_{n_2} \mathfrak{p}_s(n_2)\, \rho^{n_2}=\frac{ \mathfrak{f}_3\left(s,\frac{\rho}{12}\right)}{\mathfrak{f}_3(s)}\ ,
\label{eq:expF2}
\end{equation}
where $\rho\in ]0,1]$ is a weight per unit volume. Recall that $n_2$ may take half integer values so that the sum on the left hand side above actually
runs over all (positive) half-integers.

Let us now discuss the scaling correspondence in details\footnote{We discuss here the general case where $h\neq g$ is fixed while $g\to 1/12$.
Our arguments could be repeated {\it verbatim} to the case $h=g\to 1/12$ to explain the scaling correspondence in this case, as observed 
directly from the explicit expressions of $F(s,g,g)$ and $\mathcal{F}(S,a)$.}. Its origin is best understood by considering the all order expansion of $F(s,g,h)$ for $g\to 1/12$, namely
\begin{equation*}
F(s,g,h)=\sum_{i\geq 0} \mathfrak{f}_i(s,h)\, (1-12\, g)^{i/2}
\end{equation*}
(with $\mathfrak{f}_1(s,h)=0$ as discussed in the footnote \ref{note1}). We may indeed, via the identification $(1-12\, g)^{1/2}=(a/\sqrt{6})^2 \epsilon^2$ for $g=G(a,\epsilon)$, 
relate the scaling function $\mathcal{F}(S,a,b)$ to this expansion upon writing 
\begin{equation*}
\begin{split}
\mathcal{F}(S,a,b)&=\lim_{\epsilon\to 0}\,  \frac{1}{\epsilon^3} \, F\left(\left\lfloor S/\epsilon\right\rfloor,G(a,\epsilon),G(b,\epsilon)\right)\\
&=\lim_{\epsilon\to 0} \, \sum_{i\geq 0} \left(\frac{a}{\sqrt{6}}\right)^{2i}\, \epsilon^{2i-3}\, \mathfrak{f}_i\left(\left\lfloor S/\epsilon\right\rfloor,G(b,\epsilon)\right)\ .\\
\end{split}
\end{equation*}
Since $G(b,\epsilon)$ depends only on the product $b\, \epsilon$, the quantity $\mathfrak{f}_i\left(\left\lfloor S/\epsilon\right\rfloor,G(b,\epsilon)\right)$, which depends {\it a priori} on $S$, $b$ and $\epsilon$, is actually a function 
of the two variables $S/\epsilon$ and $b\, \epsilon$ only, 
or equivalently of the two variables  $S/\epsilon$ and $S/\epsilon\times b\, \epsilon = b\, S$. We deduce from the very existence of the scaling function above that\footnote{The fact that all the $\varphi_i$, $i\neq 1$, are not zero is verified {\it a posteriori} by the fact that $\mathcal{F}(S,a,\tau/S)$ has a 
small $S$ expansion involving non vanishing $S^{2i-3}$ coefficients for all $i\geq 0$, $i\neq 1$.}  
\begin{equation}
\mathfrak{f}_i\left(\left\lfloor S/\epsilon\right\rfloor,G(b,\epsilon)\right) \underset{\epsilon \to 0}{\sim} \left(\frac{S}{\epsilon}\right)^{2i-3} \varphi_i(b\, S)
\label{eq:fiphi}
\end{equation}
for $i\neq 1$ (while $\mathfrak{f}_1=0$) with, moreover, the direct identification
\begin{equation*}
\mathcal{F}(S,a,b)=
 \sum_{i\geq 0} \left(\frac{a}{\sqrt{6}}\right)^{2i}\, S^{2i-3}\, \varphi_i(b\, S)
\end{equation*}
with $\varphi_1=0$. 
This latter identity allows us in turn to identify the functions $\varphi_i(\tau)$ via
\begin{equation}
\varphi_i(\tau)=  \left(\frac{a}{\sqrt{6}}\right)^{-2i} [S^{2i-3}] \mathcal{F}\left(S,a,\frac{\tau}{S}\right) =  [S^{2i-3}]  \mathcal{F}\left(S,\sqrt{6},\frac{\tau}{S}\right)\ ,
\label{eq:phival}
\end{equation}
where the last term was obtained by setting $a=\sqrt{6}$ in the middle term since, $\varphi_i(\tau)$ being independent of $a$, the middle term should be too for consistency.
\vskip .3cm
We may now come back to our estimate of \eqref{eq:expF2} when $s$ is large. To obtain a
non trivial law at large $s$, we must measure $n_2$ in units of $s^4$, i.e.\ consider
the probability distribution for the rescaled volume $V$ defined by
\begin{equation*}
 V\equiv \frac{n_2}{s^4}\ .
 \end{equation*}
This law is indeed captured by choosing $\rho=e^{-\sigma/s^4}$, in which case the second argument, $\rho/12$, of the numerator
in \eqref{eq:expF2} behaves as
\begin{equation*}
\frac{1}{12}e^{-\sigma/s^4} \underset{s \to \infty }{\sim}\frac{1}{12} \left(1-\frac{\sigma}{S^4} \left(\frac{S}{s}\right)^4\right) = G\left(\frac{\sqrt{6}\, \sigma^{1/4}}{S},\frac{S}{s}\right)\ .
\end{equation*}
Taking $\epsilon=S/s$ and $b=\sqrt{6}\, \sigma^{1/4}/S$ in the above estimate \eqref{eq:fiphi} and using the identification \eqref{eq:phival}, we may now write
\begin{equation*}
\mathfrak{f}_i\left(s,\frac{1}{12}e^{-\sigma/s^4} \right) \underset{s \to \infty }{\sim} s^{2i-3}\  \varphi_i(\sqrt{6}\, \sigma^{1/4})= s^{2i-3}\ [S^{2i-3}]  \mathcal{F}\left(S,\sqrt{6},\frac{\sqrt{6}\, \sigma^{1/4}}{S}\right)\ ,
\end{equation*}
leading eventually, using \eqref{eq:largesbis} and \eqref{eq:scalF}, to
\begin{equation}
\hskip -10pt \sum_{n_2}\mathfrak{p}_s(n_2)\, e^{-\sigma\, V} \underset{s \to \infty
}{\sim} \frac{[S^3]  \mathcal{F}\left(S,\sqrt{6},\frac{\sqrt{6}\sigma^{1/4}}{S}\right)}{[S^3]  \mathcal{F}(S,\sqrt{6})}=\frac{7}{16}\,  [S^3]  \mathcal{F}\left(S,\sqrt{6},\frac{\sqrt{6}\sigma^{1/4}}{S}\right)\ .
\label{eq:exprasymp}
\end{equation}
Since we have at our disposal an explicit expression for $\mathcal{F}(S,a,b)$, this equation will give us a direct access to the desired
law for $V$. 

\subsection{Explicit expressions and plots}
\label{sec:result}
As a first, rather trivial, check of our expression \eqref{eq:exprasymp}, let us estimate the probability that $n_2$ remains finite in
infinitely large planar bi-pointed quadrangulations.
This probability is obtained by summing $\mathfrak{p}_s(n_2) $ over all allowed finite values of $n_2$, i.e.\ by setting $\rho=1$ in \eqref{eq:expF2}, i.e.\ $\sigma=0$ in \eqref{eq:exprasymp}. It therefore takes the large $s$ value
\begin{equation*}
\sum_{n_2}\mathfrak{p}_s(n_2) \underset{s \to \infty
}{\sim} \frac{ [S^3]  \mathcal{F}(S,\sqrt{6},0)}{[S^3]  \mathcal{F}(S,\sqrt{6})}=\frac{7}{16}\, [S^3]  \mathcal{F}(S,\sqrt{6},0)\ .
\end{equation*}
For $b=0$, the explicit expression for $\mathcal{F}(S,a,b)$
simplifies into
\begin{equation*}
 \mathcal{F}(S,a,0)= -\frac{36\sqrt{2}\, a^3\, \sum\limits_{m=1}^5 p_m(a\, S)\, e^{-m\, a\, S}}{(577+408\sqrt{2})\left(\sum\limits_{m=0}^2 q_m(a\, S)\, e^{-m\, a\, S}\right)^3}\ ,
\end{equation*}
where the $p_m(r)$ ($1\leq m\leq 5$) and $q_m(r)$ ($0\leq m\leq 2$) are polynomials of degree $3$ and $2$ respectively in the variable $r$, given by
\begin{equation*}
 \begin{split}
& \hskip -30pt  p_1(r)=-6 (816+577 \sqrt{2})-6 (915+647 \sqrt{2}) r-3 (618+437 \sqrt{2}) r^2-2 (99+70 \sqrt{2}) r^3\\
&\hskip -30pt   p_2(r)=-24 (222+157 \sqrt{2})-12(126+89 \sqrt{2}) r+12 (27+19 \sqrt{2}) r^2+4(24+17 \sqrt{2}) r^3\\
& \hskip -30pt  p_3(r)=-108 (4+3 \sqrt{2})-180 (3+2 \sqrt{2}) r-54 (4+3 \sqrt{2}) r^2-12(3+2\sqrt{2}) r^3\\
&\hskip -30pt   p_4(r)=-24 (-6+5 \sqrt{2})+12 (-6+\sqrt{2}) r-12 (3+\sqrt{2}) r^2-4 \sqrt{2} r^3\\
&\hskip -30pt   p_5(r)=-6(-24+17 \sqrt{2})+6(-27+19 \sqrt{2}) r-3 (-18+13 \sqrt{2}) r^2+2(-3+2 \sqrt{2}) r^3\\
&\hskip -30pt   q_0(r)=6+3 \sqrt{2} r+r^2\\
&\hskip -30pt   q_1(r)=-24(-4+3 \sqrt{2})-12 (-3+2 \sqrt{2}) r+2(-4+3 \sqrt{2}) r^2\\
&\hskip -30pt   q_2(r)=6(-17+12 \sqrt{2})-3(-24+17 \sqrt{2}) r+(-17+12 \sqrt{2}) r^2\ .\\
 \end{split}
\end{equation*}
 We have in particular the small $S$ expansion
 \begin{equation}
  \mathcal{F}(S,a,0)=\frac{1}{2 S^3}-\frac{a^4 S}{60}+\frac{a^6 S^3}{189}+O\left(S^4\right)
  \label{eq:FSazero}
 \end{equation}
 which leads to a probability that $n_2$ be finite equal to
 \begin{equation*}
\sum_{n_2}\mathfrak{p}_s(n_2) \underset{s \to \infty
}{\sim} \frac{7}{16}\, \frac{(\sqrt{6})^6}{189}=\frac{1}{2}\ .
\end{equation*}
We thus see that configurations for which the second Vorono\"\i\ cell remains finite whenever $N\to\infty$ represent at large $s$ precisely $1/2$ of
all the configurations. This is fully consistent with our result \eqref{eq:probahalf} provided that the configurations for which we found
$n_2=o(N)$ are actually configurations for which $n_2$ remains finite. Otherwise stated, configurations for which both $n_1$ and $n_2$ would
diverge at large $N$ are negligible at (large) finite $s$.

\begin{figure}
\begin{center}
\includegraphics[width=9.cm]{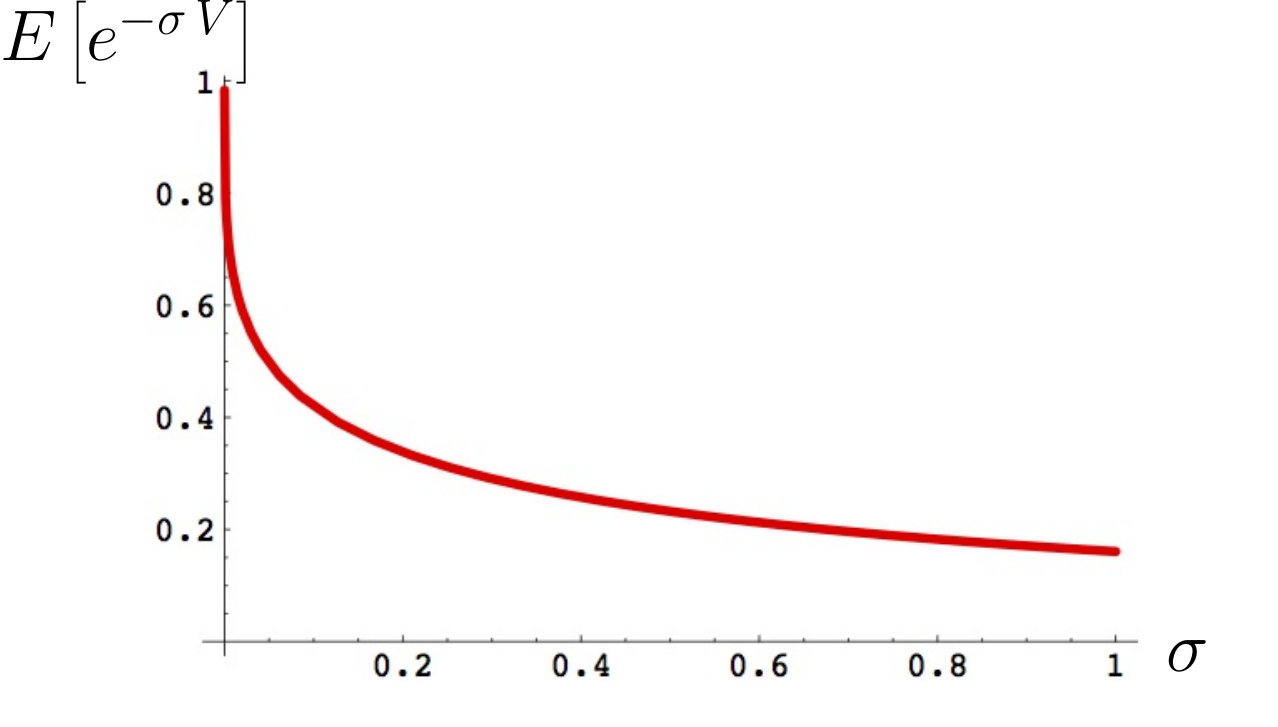}
\end{center}
\caption{A plot of the expectation value $E\left[e^{-\sigma\, V}\right]$ as given by its explicit expression \eqref{eq:Esigma}.}
\label{fig:expV}
\end{figure}
\vskip .3cm
 Beyond this first result at $\sigma=0$, we can consider, for any $\sigma \geq 0$,  the expectation value of $e^{-\sigma\, V}$ 
 for bi-pointed quadrangulations with $N\to \infty$ and finite $s$, \emph{conditioned to have their second Vorono\"\i\ cell finite}\footnote{Alternatively, we may lift this conditioning and interpret 
$E_s\left[e^{-\sigma\, V}\right]$
as the expectation value of $e^{-\sigma\, V}$ where $V$ is the rescaled volume of the smallest Vorono\"\i\ cell.
}. It is given by
 \begin{equation*}
E_s\left[e^{-\sigma\, V}\right]= \frac{ \sum\limits_{n_2}\mathfrak{p}_s(n_2)\, e^{-\sigma\, V} }{\sum\limits_{n_2}\mathfrak{p}_s(n_2)}\ ,
\end{equation*}
and has a large $s$ limiting value
\begin{equation*}
 \begin{split}
E_s\left[e^{-\sigma\, V}\right]
\underset{s \to \infty}{\sim} E\left[e^{-\sigma\, V}\right]
& =\frac{[S^3]  \mathcal{F}\left(S,\sqrt{6},\frac{\sqrt{6}\sigma^{1/4}}{S}\right)}{[S^3]  \mathcal{F}\left(S,\sqrt{6},0\right)}\\
&
= \frac{7}{8}\  [S^3]  \mathcal{F}\left(S,\sqrt{6},\frac{\sqrt{6}\sigma^{1/4}}{S}\right)\ .\\
\end{split}
\end{equation*}
From the explicit expression of $\mathcal{F}(S,a,b)$ displayed in Appendix A, we deduce after some quite heavy computation the following expression for $E\left[e^{-\sigma\, V}\right]$:
\begin{equation}
\hskip -10pt E\left[e^{-\sigma\, V}\right]=
 \frac{3}{2}\frac{P(\sigma^{1/4})+ \sum\limits_{m=1}^3\left( P_m(\sigma^{1/4},\sqrt{2})\, e^{m\sqrt{6}\, \sigma^{1/4}}
+ P_m(\sigma^{1/4},-\sqrt{2})\,e^{-m\sqrt{6}\, \sigma^{1/4}}\right)}
{\left(Q(\sigma^{1/4})\left(4+(4+3\sqrt{2})e^{\sqrt{6}\sigma^{1/4}}+(4-3\sqrt{2})
e^{-\sqrt{6}\sigma^{1/4}}\right)-12\right)^4}
\label{eq:Esigma}
\end{equation}
where $P(r)$, $Q(r)$ and $P_m(r,\gamma)$ ($m=1,2,3$) are polynomials in $r$ of degree at most $8$ (with coefficients linear in $\gamma$ for the last three 
polynomials), given explicitly by

\begin{equation*}
 \begin{split}
&\hskip -32pt  P(r)=96 \left(-252-399 \sqrt{3} r-756 r^2-161 \sqrt{3} r^3+170 r^4+153 \sqrt{3} r^5+144 r^6+22 \sqrt{3} r^7+4 r^8\right)\\
\hskip -30pt  P_1(r,\gamma)&=126 (168+85 \gamma )+63 \sqrt{3}\, r\, (867+596\gamma )+1323\, r^2\, (132+95 \gamma )+28 \sqrt{3}\, r^3\, (3153+2300 \gamma )\\
&
+24\, r^4\, (2463+1843 \gamma ) +\sqrt{3}\, r^5\, (588+905 \gamma )-36\, r^6\, (177+124 \gamma )-6 \sqrt{3}\, r^7\, (174+127 \gamma )\\
&
-36\, r^8 \, (4+3 \gamma )   
\\
\hskip -30pt   P_2(r,\gamma)&=-8 \left(63 (24+17 \gamma )+63 \sqrt{3}\, r\, (105+74 \gamma )+378\, r^2\,  (78+55 \gamma )+14\sqrt{3}\, r^3\, (1569+1108 \gamma )\right.\\
&\left.  +12 \, r^4\,  (2337+1652 \gamma )+ \sqrt{3}\, r^5\, (6954+4919 \gamma )+18\, r^6\,  (154+109 \gamma )+6 \sqrt{3} \, r^7\,  (24+17 \gamma )\right)
   \\
\hskip -30pt   P_3(r,\gamma)&=126 (24+17 \gamma )+63 \sqrt{3}\, r\, (277+196 \gamma )+189 \,r^2\, (516+365 \gamma )+28 \sqrt{3}\, r^3 \,(3399+2404 \gamma )\\
&+24 \,r^4\, (7193+5087 \gamma )+\sqrt{3}\,r^5\,  (68436+48397 \gamma )+36\, r^6\, (1465+1036 \gamma )\\&+6 \sqrt{3}\, r^7\, (1342+949 \gamma )
+12 r^8 (140+99 \gamma )\\
& \hskip -30pt  Q(r)= \left(1+\sqrt{3} r+r^2\right)\ .\\
 \end{split}
\end{equation*}
\begin{figure}
\begin{center}
\includegraphics[width=12cm]{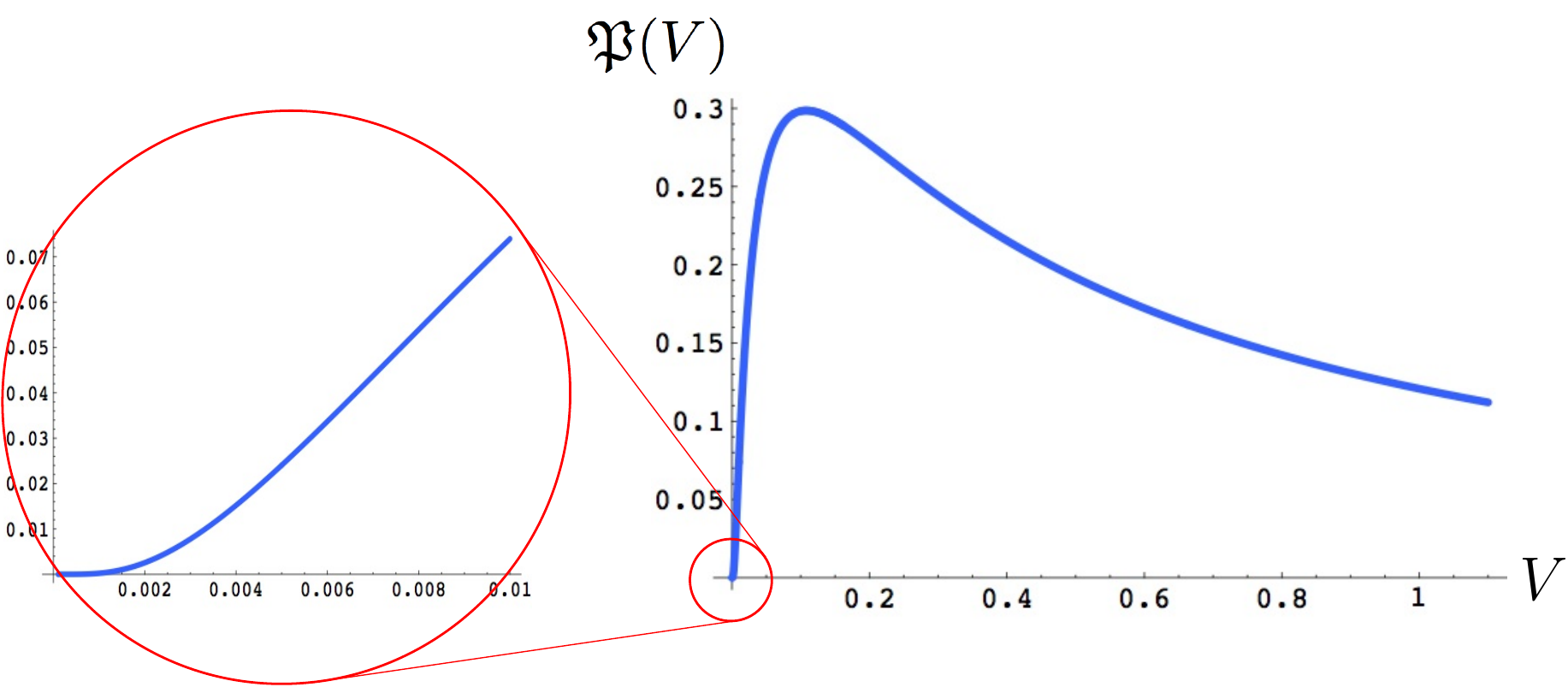}
\end{center}
\caption{A plot of the probability distribution $\mathfrak{P}(V)$ as obtained from the explicit expression \eqref{eq:Esigma} via some
numerical inverse Laplace transform. As displayed in the insert, this law becomes flat when $V\to 0$.}
\label{fig:ProbaV}
\end{figure}
The function $E\left[e^{-\sigma\, V}\right]$ is plotted in Figure \ref{fig:expV} for illustration.
\vskip .3cm
The actual probability distribution $\mathfrak{P}_s(V)$ for the rescaled volume $V$ is the inverse Laplace transform of $E_s\left[e^{-\sigma\, V}\right]$, hence its large $s$
limit $\mathfrak{P}(V)$ is given by the inverse Laplace transform of $E\left[e^{-\sigma\, V}\right]$. From the quite involved form 
\eqref{eq:Esigma} above, there is no real hope to get an explicit
expression for $\mathfrak{P}(V)$ but it may still be plotted thanks to appropriate numerical tools \cite{NILT,Ava,AVb}. The resulting shape is displayed in Figure \ref{fig:ProbaV}.
\vskip .3cm
A few analytic properties of $\mathfrak{P}(V)$ may be obtained from its explicit Laplace transform \eqref{eq:Esigma} above: in particular, we may easily find large and small $V$ 
asymptotic equivalents of $\mathfrak{P}(V)$, as discussed now.
\vskip .3cm
$\bullet$ {\bf The large volume limit.} For small $\sigma$, we have the expansion 
 \begin{figure}
\begin{center}
\includegraphics[width=8cm]{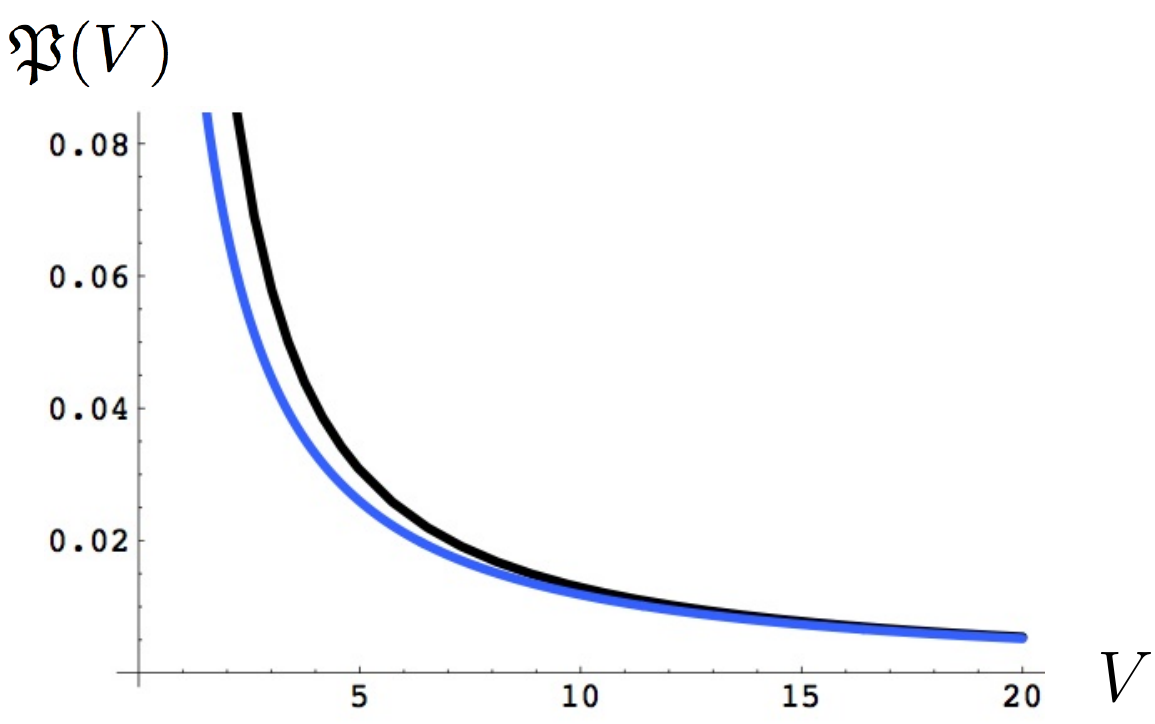}
\end{center}
\caption{Algebraic decay of $\mathfrak{P}(V)$ at large $V$: the figure displays a comparison between $\mathfrak{P}(V)$ and its large $V$
asymptotic equivalent (in black) given by \eqref{eq:PVlargeV}.}
\label{fig:ProbaVinfty}
\end{figure}
 \begin{equation*}
E\left[e^{-\sigma\, V}\right] = 1-\frac{665 \sqrt{3}}{1024} \sigma^{1/4}+\frac{49}{768\sqrt{3}}\sigma^{3/4}+\frac{63}{80}\sigma+O(\sigma^{5/4})
 \end{equation*}
 so that $E\left[e^{-\sigma\, V}\right]$ is not analytic at $\sigma=0$, with all its derivative infinite at this point. We first deduce that 
 all the (positive) moments of $\mathfrak{P}(V)$ are infinite. By a standard argument using the famous Karamata's tauberian theorem, the large $V$ tail of $\mathfrak{P}(V)$ is estimated from the leading ($\propto \sigma^{1/4}$) small $\sigma$ singularity above as
  \begin{equation}
\mathfrak{P}(V)\underset{V \to \infty}{\sim} \frac{665 \sqrt{3}}{4096 \, \Gamma(3/4)}\, \frac{1}{V^{5/4}}\ .
\label{eq:PVlargeV}
 \end{equation}
 A comparison between $\mathfrak{P}(V)$ (as obtained numerically) and its large $V$ equivalent is displayed in Figure \ref{fig:ProbaVinfty}.
\vskip .3cm 
$\bullet$ {\bf The mall volume limit.} For large $\sigma$, we have the asymptotic equivalence
 \begin{figure}
\begin{center}
\includegraphics[width=8cm]{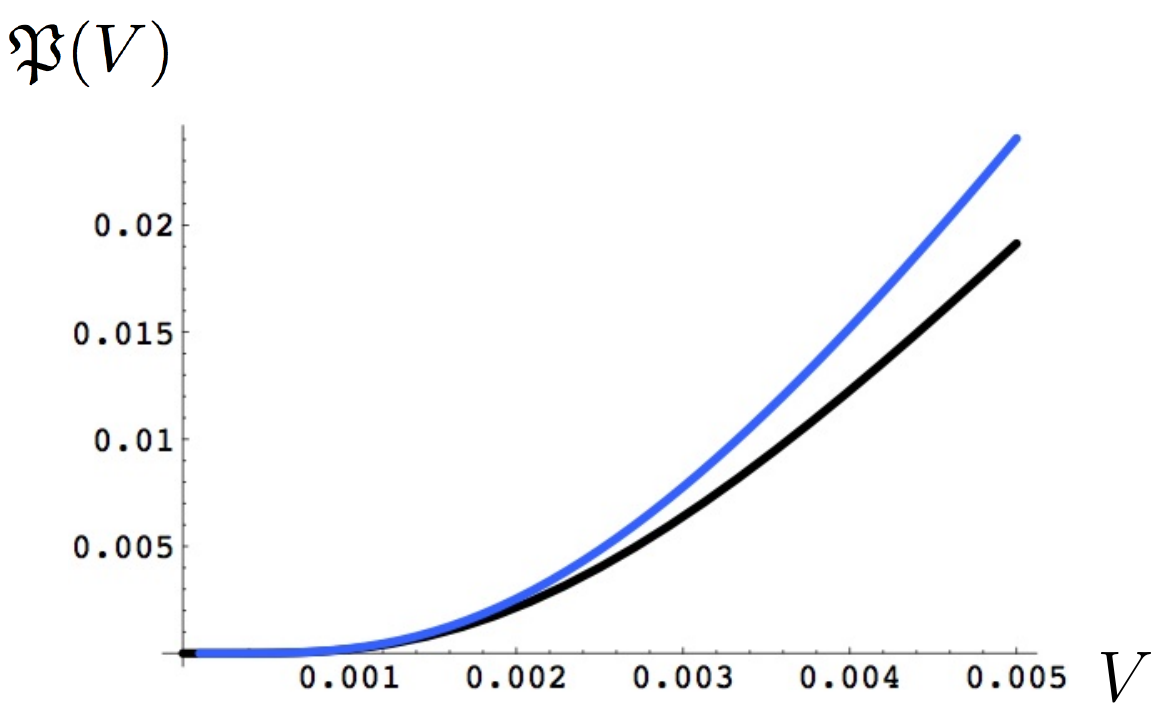}
\end{center}
\caption{The flat nature of $\mathfrak{P}(V)$ at small $V$: the figure displays a comparison between $\mathfrak{P}(V)$ and its 
small $V$ equivalent (in black) given by \eqref{eq:PVsmallV}.}
\label{fig:ProbaV0}
\end{figure}
\begin{equation*}
E\left[e^{-\sigma\, V}\right] \underset{\sigma \to \infty}{\sim} \frac{9}{2}(3\sqrt{2}-4)\, e^{-\sqrt{6}\, \sigma^{1/4}}\ .
 \end{equation*}
 By a simple saddle point calculation (see Appendix C), we deduce the small $V$ estimate
  \begin{equation}
\mathfrak{P}(V)\underset{V \to 0}{\sim} \frac{3^{11/6} (3 -2 \sqrt{2})}{2 \sqrt{\pi }}\, \frac{1}{V^{7/6}}\, e^{-\frac{3^{5/3}}{4}\frac{1}{V^{1/3}}}
\label{eq:PVsmallV}
 \end{equation} 
 which is a \emph{flat} function at $V=0$. 
 A comparison between $\mathfrak{P}(V)$ (as obtained numerically) and its small $V$ equivalent is displayed in Figure \ref{fig:ProbaV0}.
 
\section{A comparison with Vorono\"\i\ cells in infinitely large bi-pointed trees}
As an exercise, it is interesting to compare our result for $\mathfrak{P}(V)$ to that, much simpler, obtained for another family of maps, namely \emph{bi-pointed plane trees},
which are planar maps with a single face and with two marked vertices $v_1$ and $v_2$, taken again at some fixed even distance $d(v_1,v_2)=2s$
along the tree. Any such map is made of a simple path $\mathcal{P}$, formed by the edges joining $v_1$ to $v_2$, completed
by trees attached to the internal vertices of $\mathcal{P}$ on both side of the path and at its extremities $v_1$ and $v_2$. The two Vorono\"\i\ cells are now trivially
defined by splitting the tree at the ``central vertex" in $\mathcal{P}$, which is the vertex along $\mathcal{P}$ lying at distance $s$ from both $v_1$ and $v_2$ 
(there are in general two subtrees attached to this vertex and we may decide to split the tree so as to assign one of these subtrees to the first Vorono\"\i\ cell and the other subtree to the second cell). The volumes $n_1$ and $n_2$ of the two Vorono\"\i\ cells are
now measured by their number of edges and the generating function $F_{\rm tree}(s,g,h)$ enumerating these maps
with a weight $g^{n_1}h^{n_2}$ reads simply
\begin{equation*}
F_{\rm tree}(s,g,h)=\left(g\, \left(R_{\rm tree}\left(g\right)\right)^2\right)^s\ \left(h\, \left(R_{\rm tree}\left(h\right)\right)^2\right)^s
\end{equation*}
where $R_{\rm tree}(g)$ is the generating function for planted trees with a weight $g$ per edge, namely\footnote{It is solution of $R_{\rm tree}(g)=1+g\, (R_{\rm tree}(g))^2$.}
\begin{equation*}
R_{\rm tree}(g)=\frac{1-\sqrt{1-4g}}{2g}\ .
\end{equation*}
The scaling function associated with $F_{\rm tree}(s,g,h)$ is obtained by setting
\begin{equation*}
g=G_{\rm tree}(a,\epsilon)\ ,\qquad h=G_{\rm tree}(b,\epsilon)\ , \quad \hbox{where}\ \  G_{\rm tree}(c,\epsilon)\equiv \frac{1}{4}\left(1-\frac{c^2}{4}\epsilon^2\right)\ ,
\end{equation*}
and reads simply
\begin{equation*}
\mathcal{F}_{\rm tree}(S,a,b)=\lim_{\epsilon\to 0}\,  F_{\rm tree}\left(\left\lfloor S/\epsilon\right\rfloor,G_{\rm tree}(a,\epsilon),G_{\rm tree}(b,\epsilon)\right)
=e^{-(a+b)S}\ .
\end{equation*}
By repeating and adapting the arguments of previous sections, here with leading singularities of type $(1-4g)^{1/2}$, we can find the large $N$ ($=$ total number of edges)
asymptotic law  $\mathfrak{P}_{\rm tree}(V_{\rm tree})$ for the rescaled volume
\begin{equation*}
 V_{\rm tree}\equiv\frac{n_2}{s^2}
 \end{equation*}
among bi-pointed trees with $d(v_1,v_2)=2s$, conditioned to have their second Vorono\"\i\ cell finite (which again represent $1/2$ of all bi-pointed trees with fixed $s$). For large $s$, we find (with obvious notations) the expectation value
\begin{equation*}
 E_{\rm tree}\left[e^{-\sigma\, V_{\rm tree}}\right]=\frac{[S]  \mathcal{F}_{\rm tree}\left(S,2,\frac{2\sigma^{1/2}}{S}\right)}{[S]  \mathcal{F}_{\rm tree}\left(S,2,0\right)}
 =e^{-2\sigma^{1/2}}\ .
 \end{equation*}
 We may now deduce by inverse Laplace transform the exact law for $V_{\rm tree}$
 \begin{equation*}
 \mathfrak{P}_{\rm tree}(V_{\rm tree})=\frac{1}{\sqrt{\pi}\,  V_{\rm tree}^{3/2}}e^{-\frac{1}{V_{\rm tree}}}\ ,
 \end{equation*}
 which is nothing but a simple L\'evy distribution. In particular, all its (positive) moments are infinite. 

This distribution is a particular member of the more general family of \emph{one-sided L\'evy distributions with parameter $\alpha$},
namely distributions whose Laplace transform is $e^{-{\rm const.}\, p^\alpha}$. For $0<\alpha<1$, such distributions
are \emph{flat} at small volume $V$ and vanish as 
 \begin{equation*}
 \frac{\hbox{const.}}{V^{(2-\alpha)/(2(1-\alpha))}}\, {\rm exp}\left(-\frac{\hbox{const.'}}{V^{\alpha/(1-\alpha)}}\right)\ .
 \end{equation*}
For large volume, they present a fat tail with an algebraic decay of the form $1/V^{1+\alpha}$.
 The simple law $ \mathfrak{P}_{\rm tree}(V)$ for trees corresponds precisely to the situation where $\alpha=1/2$. 

As for the distribution $\mathfrak{P}(V)$ of quadrangulations, it is obviously not a L\'evy distribution but 
its small and large $V$ behaviors are nevertheless similar to those obtained for a L\'evy distribution with $\alpha=1/4$. 

Clearly,  the value of $\alpha$ appearing in both the small and large $V$ asymptotics is related to the fractal dimension $D$ of the maps 
at hand ($D=2$ for trees and $D=4$ for quadrangulations) via
\begin{equation*}
\alpha=\frac{1}{D}\ .
\end{equation*}
It is tempting to conjecture that the above forms for small and large $V$ asymptotics should be generic and hold for other families of maps, 
 possibly within more involved universality classes with more general fractal dimensions, hence more general values of $\alpha$.  
 
\section{Asymmetric Vorono\"\i\ cells}

In Section \ref{sec:result}, we estimated the large $N$ asymptotic proportion of bi-pointed planar quadrangulations for
which the second Vorono\"\i\ cell is finite. The obtained value $1/2$ is trivial by symmetry if we assume that configurations
for which both Vorono\"\i\ cells become infinite are negligible (this latter property being {\it de facto} proven by the result itself).
Note that, in this respect, our computation was performed here in the ``worst" situation where the value of the distance $2s$ between the two marked 
vertices is large. 

We may now explicitly break the symmetry and define \emph{asymmetric Vorono\"\i\ cells} upon introducing some bias 
in the measurement of distances. The bijection between bi-pointed planar quadrangulations and planar i-l.2.f.m is indeed only one
particular instance of the Miermont bijection. The Miermont bijection allows us to introduce more generally what are called \emph{delays}, which are integers associated with the marked vertices and allow for some asymmetry in the evaluation of distances \cite{Miermont2009}. In the case of two marked vertices, two delays may in principle be introduced but, in practice, only their difference (called $\theta$ below) does matter.
In the presence of delays, the resulting image of the bi-pointed quadrangulation is again a two-face map, but now with
a more general labelling of its vertices by integers. If we insist on keeping a distance $d(v_1,v_2)=2s$ between the marked vertices
in the original quadrangulation, the labelling of the two-face map, which now involves some additional integer parameter $\theta$, is characterized by the following 
four properties:  
\begin{enumerate}[$(\hbox{L}_1)$]
\item{labels on adjacent vertices differ by $0$ or $\pm1$;}
\item{the minimum label for the set of vertices incident to $f_1$ is $1-\theta$;}
\item{the minimum label for the set of vertices incident to $f_2$ is $1+\theta$;}
\item{the minimum label for the set of vertices incident to $\mathcal{L}$ is $s$;}
\end{enumerate}
if $f_1$ and $f_2$ denote the two faces of the map and $\mathcal{L}$ the loop made of edges incident to both faces. 
The Miermont bijection is a one-to-one correspondence between bi-pointed planar quadrangulations with $d(v_1,v_2)=2s$ 
and planar two-face maps with a labelling satisfying $(\hbox{L}_1)-(\hbox{L}_4)$ above for any \emph{fixed} $\theta$ 
in the range \cite{Miermont2009}
\begin{equation*}
-s<\theta<s\ .
\end{equation*}
All the vertices $v$ of the original quadrangulation but $v_1$ and $v_2$ are recovered in the two-face map, and 
their label is related to the distance in the quadrangulation via
\begin{equation*}
\ell(v)=\min(d(v,v_1)-\theta,d(v,v_2)+\theta)\ .
\end{equation*}
Again the two domains of the original quadrangulation covered by $f_1$ and $f_2$ respectively
(upon drawing the quadrangulation and its image via the bijection on top of each other) 
naturally define two cells in the map. For some generic $\theta$, those are however 
\emph{asymmetric Vorono\"\i\ cells} with the following properties: the first cell now contains all the vertices  $v$ 
such that 
\begin{equation*}
d(v_,v_1) < d(v,v_2)+ 2\theta \quad (\hbox{cell $1$})
\end{equation*}
(this includes the vertex $v_1$), as well as a number of vertices satisfying $d(v_,v_1)=d(v,v_2)+2\theta$. The second cell contains all the vertices $v$  such that
\begin{equation*}
d(v_,v_1) > d(v,v_2)+ 2\theta \quad (\hbox{cell $2$})
\end{equation*}
(including $v_2$) as well as a number of vertices satisfying $d(v_,v_1)=d(v,v_2)+2\theta$. In particular, the loop $\mathcal{L}$, whose vertices
belong to both $f_1$ and $f_2$, contains only
vertices satisfying  $d(v_,v_1)=d(v,v_2)+2\theta$.

Taking $\theta> 0$ therefore ``favors" cell $1$ whose volume is, on average, larger than that of cell $2$.
A control on these volumes is again obtained directly via the bijection by assigning a weight $g$ per edge in $f_1$, $h$ per edge in $f_2$
and $\sqrt{g h}$ per edge along $\mathcal{L}$. The corresponding generating function reads then 
\begin{equation*}
\hskip -1.cm F(s,\theta,g,h)= \Delta_u\Delta_v \log(X_{u,v}(g,h))\Big\vert_{u=s+\theta \atop v=s-\theta} =\log\left(\frac{X_{s+\theta,s-\theta}(g,h)X_{s+\theta-1,s-\theta-1}(g,h)}{X_{s+\theta-1,s-\theta}(g,h)X_{s+\theta,s-\theta-1}(g,h)}\right)\ ,
\end{equation*}
giving rise to a scaling function $\mathcal{F}(S,\Theta,a,b)$ via
\begin{equation*}
\begin{split}
\hskip -1.cm F(\left\lfloor S/\epsilon\right\rfloor,\left\lfloor \Theta/\epsilon\right\rfloor,G(a,\epsilon),G(b,\epsilon)) &  = \epsilon^3\, \mathcal{F}(S,\Theta,a,b) +O(\epsilon^4)\\
\hbox{where} & \quad \mathcal{F}(S,\Theta,a,b) = \frac{1}{3}\ \partial_U\partial_V x(U,V,a,b)\Big\vert_{U=S+\Theta\atop V=S-\Theta}\ .\\
\end{split}
\end{equation*}
Defining the \emph{asymmetry factor} $\omega$ by 
\begin{equation*}
\omega\equiv \frac{\theta}{s}\ , \quad -1\leq \omega \leq 1\ ,
\end{equation*}
the local limit of configurations with fixed $s$ and $\omega$ is, at large $s$, encoded in the small $S$ expansion
of  $\mathcal{F}(S,\omega\, S, a,b)$. In particular, 
we easily get the $b=0$ small $S$ expansion generalizing \eqref{eq:FSazero}
 \begin{equation*}
 \begin{split}
  \mathcal{F}(S,\omega\, S, a,0)&=\frac{1}{2 S^3}-\frac{a^4\, S}{480} (1+\omega )^3 \left(8-9 \omega +3 \omega ^2\right)\\
  &+\frac{a^6\, S^3}{6048}(1+\omega )^3 \left(32-33 \omega +3 \omega ^2+9 \omega ^3-3 \omega ^4\right)+O[S]^4\ . \\
  \end{split}
 \end{equation*}
By the same argument as in Section \ref{sec:result}, we directly read from the $S^3$ coefficient of this expansion the large $s$ probability $\Pi(\omega)$ that,
for $N\to \infty$, the volume of the second (now asymmetric with a fixed value of the asymmetry factor $\omega$) Vorono\"\i\ cell remains finite:
\begin{equation}
\Pi(\omega)=\frac{1}{64}\, (1+\omega )^3 \left(32-33 \omega +3 \omega ^2+9 \omega ^3-3 \omega ^4\right)\ .
\label{eq:Pialpha}
\end{equation}
 \begin{figure}
\begin{center}
\includegraphics[width=7cm]{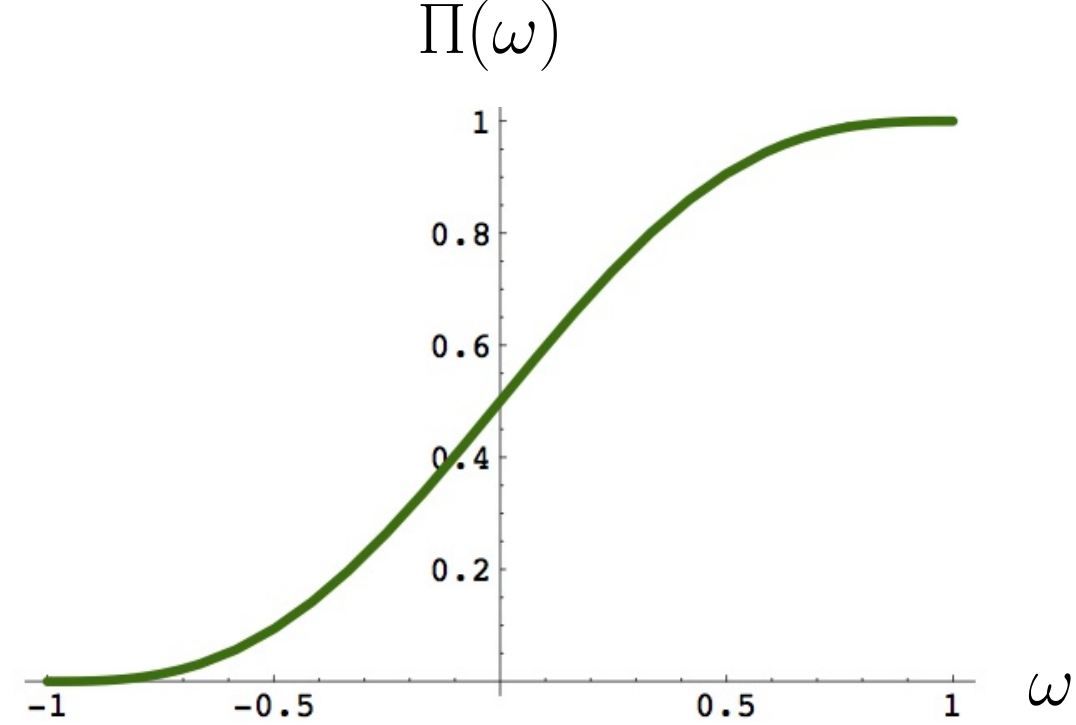}
\end{center}
\caption{A plot of the probability $\Pi(\omega)$ that the volume of the second Vorono\"\i\ cell remains finite when $N\to \infty$,
as a function of the asymmetry factor $\omega$, as given by \eqref{eq:Pialpha}.}
\label{fig:Probalpha}
\end{figure}
This probability is displayed in Figure \ref{fig:Probalpha}. It satisfies of course $\Pi(\omega)+\Pi(-\omega)=1$ as expected by symmetry, and
in particular $\Pi(0)=1/2$ for the symmetric case. For, say $\omega\to 1$, the second Vorono\"\i\ cell, which is maximally unfavored by the asymmetry
is finite with probability $\Pi(1)=1$.

For a better understanding of the meaning of $\Pi(\omega)$, we may quote Miermont in \cite{Miermont2009} and {\it ``let water flow at unit speed from the sources}" $v_1$ and $v_2$ at given mutual distance $2s$ {\it ``in such a way that the water starts diffusing from}" $v_1$ 
at time $-\omega\, s$, from $v_2$ at time $+\omega\, s$, {\it ``and takes unit time to go through an edge. When water currents emanating from different edges meet at a vertex (whenever the water initially comes from the same source of from different sources),
they can go on flowing into unvisited edges only\footnote{Miermont's bijection is so designed that meeting currents {\it ``can go on flowing into unvisited edges only respecting the rules of a roundabout, i.e.
edges that can be attained by turning around the vertex counterclockwise and not crossing any
other current."}} (...). The process ends when the water cannot flow any more (...)."}  In this language, Vorono\"\i\ cells correspond to domains covered by a given current. When the map volume tends to infinity, only one of the currents flows all the way to infinity, the other current
remaining trapped within a finite region. We may then view $\Pi(\omega)$ as
the probability for the second current (emanating from the source $v_2$) to remain trapped within a finite domain or, equivalently, as the probability
for the first current (emanating from the source $v_1$) to escape to infinity. 

\section{conclusion} 

In this paper, we computed explicitly the value of the expectation value $E\left[e^{-\sigma\, V}\right]$ for infinitely large bi-pointed planar maps, where
$V$ is the rescaled volume of that of the two Vorono\"\i\ cells which remains finite. This law describes Vorono\"\i\ cells constructed from
randomly picked vertices at a prescribed finite mutual distance, and in the limit where this distance is large. Although it may look quite involved,
the expression \eqref{eq:Esigma} is nevertheless expected to be \emph{universal} since its derivation
entirely relies on properties of scaling functions, which are in fact characteristic of the Brownian map rather than the specific realization at hand
(here quadrangulations). 
In other words, we expect that the same expression \eqref{eq:Esigma}, \emph{up to
some possible non-universal normalization} for the parameter $\sigma$, would be obtained, in the same regime, 
for all bi-pointed planar map families in the universality class of 
so-called pure gravity\footnote{This includes maps with \emph{bounded} face degrees with possible degree dependent weights, as well as 
maps with unbounded face degrees with degree dependent weights which restrain the proliferation of large faces.}.

We then deduced from this result a number of features of the associated, universal, volume probability distribution 
$\mathfrak{P}(V)$, such as its large and small $V$ behaviors \eqref{eq:PVlargeV} and \eqref{eq:PVsmallV}. 
Although we have not been able to give a tractable explicit formula for this law for arbitrary $V$, we thereby showed that
its nature is comparable to that of a simple one-sided L\'evy distribution with parameter $\alpha=1/4$.
\vskip .3cm
Let us conclude by briefly discussing Vorono\"\i\ cells, now in the so called \emph{scaling regime}: here we continue to fix the mutual
distance $2s$ between the marked vertices but we now let $s$ and $N$ tend simultaneously to infinity with the ratio $S=s/N^{1/4}$
kept finite. In this regime, the fraction $\phi=n_2/N$ of the total volume spanned by the second Vorono\"\i\ cell is again a good measure
of the cell volume distribution and the asymptotic law\footnote{We use a slightly different notation $\mathcal{P}_{\{S\}}(\phi)$ with curly brackets 
to distinguish this law from the local limit law $\mathcal{P}_s(\phi)$ at fixed $s$.} $\mathcal{P}_{\{S\}}(\phi)$ for $\phi$ at fixed $S$ may be obtained as in Section 3.1 via
\begin{equation*}
\int_{0}^{1}d\phi\, \mathcal{P}_{\{S \}}(\phi)\, e^{\mu \, \phi}= \lim_{N\to \infty} \frac{[g^N]F\left(\left\lfloor S N^{1/4}\right\rfloor,g,g\, e^{\mu/N}\right)}{[g^N]F\left(\left\lfloor S N^{1/4}\right\rfloor,g,g\right)}\ .
\end{equation*}
As explained in Appendix B in the context of the local limit law $\mathcal{P}_s(\phi)$, the $g^N$ coefficient of the numerator above may be extracted via some simple contour integral which, upon taking $g=G(a,1/N^{1/4})$, involves at large $N$ the scaling function $\mathcal{F}(S,a,b)$ at  $b=(a^4-36\mu)^{1/4}$. This leads to 
\begin{equation*}
\int_{0}^{1}d\phi\, \mathcal{P}_{\{S \}}(\phi)\, e^{\mu \, \phi}=\frac{\int_{\mathcal{C}_\mu} da\,  a^3
\, e^{a^4/36} \mathcal{F}(S,a,(a^4-36 \mu)^{1/4})
}{\int_{\mathcal{C}_0} da\,  a^3
\, e^{a^4/36} \mathcal{F}(S,a,a)}\ ,
\end{equation*}
where the integral over $a$ is over some appropriate contour $\mathcal{C}_\mu$ depending on $\mu$ (see Appendix B). 
Due to the involved expression of $\mathcal{F}(S,a,b)$, we were not able to perform the above contour integrals 
for general $S$. Still, for $S\to 0$, we recover precisely via a small $S$ expansion the result of Section 3.2 (as expected from the  scaling correspondence), now in the form
 \begin{equation*}
\lim_{S\to 0} \int_{0}^{1}d\phi\, \mathcal{P}_{\{S \}}(\phi)\, e^{\mu \, \phi}= \frac{1}{2} (1+e^\mu)
\quad \Leftrightarrow  \quad 
\lim_{S\to 0} \mathcal{P}_{\{S \}}(\phi)= \frac{1}{2} (\delta(\phi)+\delta(\phi-1))\ .
\end{equation*}
As for the $S\to\infty$ limit, it may be obtained as follows (here we simply give a sketch of the calculation and leave to the reader
the task of filling the gaps): from its explicit expression, we have
\begin{equation*}
\mathcal{F}(S,a,b)\underset{S \to \infty}{\sim} \mathfrak{Z}(a,b) e^{-(a+b)S}\ ,
\end{equation*}
where the value of the coefficient $\mathfrak{Z}(a,b)$ may easily be obtained but is unimportant for our calculation
(apart from the, easily verified property that it has a non-zero limit $\mathfrak{Z}(a,a)$ when $b\to a$).
For large $S$, the above contour integrals may be evaluated by a saddle-point method.  For the denominator (with $b=a$),
we write
\begin{equation*}
\frac{d}{da}\left(\frac{a^4}{36}-2a\, S\right)\Big\vert_{a=a^*}=0\quad \Leftrightarrow \quad (a^*)^4=(18 S)^{4/3}
\end{equation*}
and deform the contour $\mathcal{C}_0$ so as to pass via the positive real saddle point at $a^*=(18 S)^{1/3}$.
This gives a denominator\footnote{Since the denominator is directly proportional to the $S$-dependent two-point function (i.e.\ the distance profile)
in large maps, we recognize here the well-known Fisher's law which states that this function decays at large $S$ 
like $e^{-\hbox{const.} S^\delta}$ with the exponent $\delta=1/(1-1/D)=4/3$ at $D=4$.} proportional at large $S$ to  $\mathfrak{Z}(a^*,a^*)e^{-(18 S)^{4/3}/12}$.
As for the numerator, it is dominated by the same saddle point but the replacement
\begin{equation*}
b\to \left((a^*)^4-36\mu\right)^{1/4}= (18 S)^{1/3}-\frac{\mu}{2S}+O\left(\frac{1}{S^2}\right)
\end{equation*}
creates a $\mu$-dependent correction. At large $S$, this leads to a numerator now proportional to  $\mathfrak{Z}(a^*,a^*)e^{-(18 S)^{4/3}/12+\mu/2}$ with the same 
proportionality constant as for the denominator, giving eventually
 \begin{equation*}
\lim_{S\to \infty} \int_{0}^{1}d\phi\, \mathcal{P}_{\{S \}}(\phi)\, e^{\mu \, \phi}=e^{\mu/2}
\quad \Leftrightarrow  \quad 
\lim_{S\to \infty} \mathcal{P}_{\{S \}}(\phi)= \delta(\phi-1/2)\ .
\end{equation*}
This result is quite natural since, heuristically, the limit $S\to \infty$ describes maps with an \emph{elongated shape}, with the
two marked vertices sitting at its extremities. The frontier between the Vorono\"\i\ cells
for such an elongated map is typically a small cycle sitting halfway along the elongated direction, hence splitting the map into
two domains of the same volume, equal to half the total volume. It would be interesting to visualize and follow the continuous passage, for increasing
$S$, of the distribution $\mathcal{P}_{\{S \}}(\phi)$ from its $S\to 0$ to its $S\to \infty$ limit above and to better understand how its average
over arbitrary $S$ (properly weighted by the $S$-dependent two-point function, i.e.\ the distance profile of the Brownian map) creates the uniform distribution for $\phi\in[0,1]$.

\section*{Acknowledgements} 
The author acknowledges the support of the grant ANR-14-CE25-0014 (ANR GRAAL).

\appendix
\section{Expression for the scaling function $\mathcal{F}(S,a,b)$} 
The scaling function $x(S,T,a,b)$ was computed in \cite{G17a} as the appropriate solution of \eqref{eq:eqforxstab}. Its explicit expression is quite heavy and is not reproduced here. From this expression, we may obtain $\mathcal{F}(S,a,b)$ directly via \eqref{eq:xtoF}. It takes the following form:
\begin{equation*}
\mathcal{F}(S,a,b)=-\frac{e^{-(a+b)\, S}}{6} \  \frac{\mathfrak{T}(e^{-a\, S},e^{-b\, S},a,b)}{\mathfrak{D}(a,b)}\,\left( \frac{\mathfrak{E}(a,b)}{\mathfrak{U}(e^{-a\, S},e^{-b\, S},a,b)
}\right)^3\ ,
\end{equation*}
where, introducing the notation
\begin{equation*}
c \equiv \sqrt{\frac{a^2+b^2}{2}}\ ,
\end{equation*}
we have explicitly
\begin{equation*}
\begin{split}
\hskip -1.2cm \mathfrak{E}(a,b)&= 6\, a\, b\,  (a\!-\!b)^2  (a\!+\!b) \left(2 a^2\!+\!b^2\right)\left(a^2\!+\!2 b^2\right)\ ,\\
\hskip -1.2cm \mathfrak{D}(a,b)&=(a\!+\!2 c) (b\!+\!2 c) \left(5 a^3\!+\!7 a^2 c\!+\! 4 a b^2 \!+\!2 b^2 c\right) \left(4a^2 b\!+\!2a^2 c\!+\!5 b^3 \!+\!7 b^2 c\right)\\
&\times \left(17 a^2 (a^2 \!+\! b^2)\!+\!12 a \left(2 a^2\!+\!b^2\right) c \!+\!2 b^4\right)
   \left(2 a^4\!+\!12 b \left(a^2\!+\!2 b^2\right) c\!+\! 17 b^2 (a^2 \!+\! b^2)\right)\ ,
\\
\end{split}
\end{equation*}
while  $\mathfrak{T}(\sigma,\tau,a,b)$ and $\mathfrak{U}(\sigma,\tau,a,b)$ are polynomials of respective degree $4$ and $2$ in both $\sigma$ and $\tau$,
namely
\begin{equation*}
\mathfrak{T}(\sigma,\tau,a,b)=\sum_{i=0}^{4}\sum_{j=0}^4 t_{i,j}\, \sigma^i\tau^j  \quad \hbox{and}\quad
\mathfrak{U}(\sigma,\tau,a,b) =\sum_{i=0}^{2}\sum_{j=0}^2 u_{i,j}\, \sigma^i\tau^j \ .
\end{equation*}
The coefficients $t_{i,j}\equiv t_{i,j}(a,b)$ and $u_{i,j}\equiv u_{i,j}(a,b)$ may be written for convenience as sums of two contributions:
\begin{equation*}
t_{i,j}=t_{i,j}^{(0)}+c\, t_{i,j}^{(1)}\ , \qquad
u_{i,j}=u_{i,j}^{(0)}+c\, u_{i,j}^{(1)}\ ,
\end{equation*}
where we have the explicit expressions
\begin{equation*}
\begin{split}
& \hskip -1.2cm t_{0,0}^{(0)}=(a\!+\!b)^3 (396 a^{10}\!+\!1448 a^9 b\!+\!3672 a^8 b^2\!+\!6520 a^7 b^3\!+\!9135 a^6 b^4\!+\!10146 a^5 b^5\!+\!9135 a^4 b^6\\
&\!+\!6520 a^3 b^7\!+\!3672 a^2 b^8\!+\!1448 a b^9\!+\!396 b^{10})\\ 
& \hskip -1.2cm t_{0,1}^{(0)}=-4 (a^2\!-\!b^2)^2 (198 a^9\!+\!502 a^8 b\!+\!1099 a^7 b^2\!+\!1551 a^6 b^3\!+\!1806 a^5 b^4\!+\!1596 a^4 b^5\\ &
\!+\!1128 a^3 b^6\!+\!596 a^2 b^7\!+\!224 a b^8\!+\!48 b^9)\\ 
& \hskip -1.2cm t_{0,2}^{(0)}=-6 b (a^2\!-\!b^2) \left(2 a^2\!+\!b^2\right) (198 a^8\!+\!280 a^7 b\!+\!611 a^6 b^2\!+\!584 a^5 b^3\!+\!599 a^4 b^4\!+\!368 a^3 b^5\\&
\!+\!200 a^2 b^6\!+\!64 a b^7\!+\!12 b^8)\\ 
& \hskip -1.2cm t_{0,3}^{(0)}=4 a (a^2\!-\!b^2)^2  \left(2 a^2\!+\!b^2\right) \left(99 a^6\!+\!29 a^5 b\!+\!186 a^4 b^2\!+\!40 a^3 b^3\!+\!104 a^2 b^4\!+\!12 a b^5\!+\!16 b^6\right)\\ 
& \hskip -1.2cm t_{0,4}^{(0)}=-(a\!-\!b)^3 \left(2 a^2\!+\!b^2\right)^2 \left(99 a^6-82 a^5 b\!+\!191 a^4 b^2-120 a^3 b^3\!+\!104 a^2 b^4-40 a b^5\!+\!12 b^6\right)\\ 
\end{split}
\end{equation*}
\begin{equation*}
\begin{split}
& \hskip -1.2cm t_{1,0}^{(0)}= -4 (a^2\!-\!b^2)^2(48 a^9\!+\!224 a^8 b\!+\!596 a^7 b^2\!+\!1128 a^6 b^3\!+\!1596 a^5 b^4\!+\!1806 a^4 b^5\\ &\!+\!1551 a^3 b^6\!+\!1099 a^2 b^7\!+\!502 a b^8\!+\!198 b^9)\\ 
& \hskip -1.2cm t_{1,1}^{(0)}=8 (a\!+\!b)^3 \left(a^4\!+\!7 a^2 b^2\!+\!b^4\right) \left(48 a^6\!+\!74 a^5 b\!+\!168 a^4 b^2\!+\!149 a^3 b^3\!+\!168 a^2 b^4\!+\!74 a b^5\!+\!48 b^6\right) \\ 
& \hskip -1.2cm t_{1,2}^{(0)}=-24 b (a^2\!-\!b^2)^2 \left(2 a^2\!+\!b^2\right) \left(24 a^6\!+\!34 a^5 b\!+\!62 a^4 b^2\!+\!54 a^3 b^3\!+\!43 a^2 b^4\!+\!20 a b^5\!+\!6 b^6\right)\\ 
& \hskip -1.2cm t_{1,3}^{(0)}=-8 a (a\!-\!b)^2 \left(2 a^2\!+\!b^2\right) \left(a^4\!+\!7 a^2 b^2\!+\!b^4\right) \left(24 a^4\!+\!7 a^3 b\!+\!33 a^2 b^2\!+\!6 a b^3\!+\!10 b^4\right) \\ 
& \hskip -1.2cm t_{1,4}^{(0)}= 4 (a^2\!-\!b^2)^2 \left(2 a^2\!+\!b^2\right)^2 \left(12 a^5-22 a^4 b\!+\!27 a^3 b^2-27 a^2 b^3\!+\!14 a b^4-6 b^5\right)\\
\end{split}
\end{equation*}
\begin{equation*}
\begin{split} 
& \hskip -1.2cm t_{2,0}^{(0)}=6 a (a^2\!-\!b^2) \left(a^2\!+\!2 b^2\right) (12 a^8\!+\!64 a^7 b\!+\!200 a^6 b^2\!+\!368 a^5 b^3\!+\!599 a^4 b^4\!+\!584 a^3 b^5\\
&\!+\!611 a^2 b^6\!+\!280 a b^7\!+\!198 b^8) \\ 
& \hskip -1.2cm t_{2,1}^{(0)}= -24 a (a^2\!-\!b^2)^2  \left(a^2\!+\!2 b^2\right) \left(6 a^6\!+\!20 a^5 b\!+\!43 a^4 b^2\!+\!54 a^3 b^3\!+\!62 a^2 b^4\!+\!34 a b^5\!+\!24 b^6\right)\\ 
& \hskip -1.2cm t_{2,2}^{(0)}= 0\\ 
& \hskip -1.2cm t_{2,3}^{(0)}=24 a^2 (a^2\!-\!b^2)^2 \left(2 a^2\!+\!b^2\right) \left(a^2\!+\!2 b^2\right) \left(3 a^3-2 a^2 b\!+\!2 a b^2-2 b^3\right)\\ 
& \hskip -1.2cm t_{2,4}^{(0)}=-6 a (a^2\!-\!b^2)  \left(2 a^2\!+\!b^2\right)^2 \left(a^2\!+\!2 b^2\right) \left(3 a^4-8 a^3 b\!+\!11 a^2 b^2-8 a b^3\!+\!6 b^4\right)\\
\end{split}
\end{equation*}
\begin{equation*}
\begin{split} 
& \hskip -1.2cm t_{3,0}^{(0)}=4 b (a^2\!-\!b^2)^2  \left(a^2\!+\!2 b^2\right) \left(16 a^6\!+\!12 a^5 b\!+\!104 a^4 b^2\!+\!40 a^3 b^3\!+\!186 a^2 b^4\!+\!29 a b^5\!+\!99 b^6\right)\\ 
& \hskip -1.2cm t_{3,1}^{(0)}= -8 b (a\!-\!b)^2 \left(a^2\!+\!2 b^2\right) \left(a^4\!+\!7 a^2 b^2\!+\!b^4\right) \left(10 a^4\!+\!6 a^3 b\!+\!33 a^2 b^2\!+\!7 a b^3\!+\!24 b^4\right)\\ 
& \hskip -1.2cm t_{3,2}^{(0)}=-24 b^2 (a^2\!-\!b^2)^2\left(2 a^2\!+\!b^2\right) \left(a^2\!+\!2 b^2\right) \left(2 a^3\!-\!2 a^2 b\!+\!2 a b^2\!-\!3 b^3\right)\\ 
& \hskip -1.2cm t_{3,3}^{(0)}= -8 a b (a\!+\!b)^3 \left(2 a^2\!+\!b^2\right) \left(a^2\!+\!2 b^2\right) \left(a^4\!+\!7 a^2 b^2\!+\!b^4\right)\\ 
& \hskip -1.2cm t_{3,4}^{(0)}= 4 b (a^2\!-\!b^2)^2  \left(2 a^2\!+\!b^2\right)^2 \left(a^2\!+\!2 b^2\right) \left(2 a^2-a b\!+\!3 b^2\right)\\
\end{split}
\end{equation*}
\begin{equation*}
\begin{split} 
& \hskip -1.2cm t_{4,0}^{(0)}=(a\!-\!b)^3 \left(a^2\!+\!2 b^2\right)^2 \left(12 a^6-40 a^5 b\!+\!104 a^4 b^2-120 a^3 b^3\!+\!191 a^2 b^4-82 a b^5\!+\!99 b^6\right)\\ 
& \hskip -1.2cm t_{4,1}^{(0)}=-4 (a^2-b^2)^2  \left(a^2+2 b^2\right)^2 \left(6 a^5-14 a^4 b+27 a^3 b^2-27 a^2 b^3+22 a b^4-12 b^5\right)\\ 
& \hskip -1.2cm t_{4,2}^{(0)}= 6 b(a^2\!-\!b^2)  \left(2 a^2\!+\!b^2\right) \left(a^2\!+\!2 b^2\right)^2 \left(6 a^4-8 a^3 b\!+\!11 a^2 b^2-8 a b^3\!+\!3 b^4\right)\\ 
& \hskip -1.2cm t_{4,3}^{(0)}= 4 a (a^2-b^2)^2 \left(2 a^2+b^2\right) \left(a^2+2 b^2\right)^2 \left(3 a^2-a b+2 b^2\right)\\ 
& \hskip -1.2cm t_{4,4}^{(0)}= -(a\!+\!b)^3 \left(2 a^2\!+\!b^2\right)^2 \left(a^2\!+\!2 b^2\right)^2 \left(3 a^2\!+\!2 a b\!+\!3 b^2\right)\\ 
\end{split}
\end{equation*}
and
\begin{equation*}
\begin{split}
& \hskip -1.2cm t_{0,0}^{(1)}=4 (a\!+\!b)^4 \left(10 a^4\!+\!18 a^3 b\!+\!25 a^2 b^2\!+\!18 a b^3\!+\!10 b^4\right) \left(14 a^4\!+\!12 a^3 b\!+\!29 a^2 b^2\!+\!12 a b^3\!+\!14 b^4\right)\\ 
& \hskip -1.2cm t_{0,1}^{(1)}=-4 (a^2\!-\!b^2)^2(280 a^8\!+\!710 a^7 b\!+\!1414 a^6 b^2\!+\!1839 a^5 b^3\!+\!1881 a^4 b^4\!+\!1428 a^3 b^5\!+\!812 a^2 b^6\\&\!+\!316 a b^7\!+\!68 b^8)\\ 
& \hskip -1.2cm t_{0,2}^{(1)}=-24 b (a^2\!-\!b^2)\left(2 a^2\!+\!b^2\right) \left(10 a^3\!+\!7 a^2 b\!+\!8 a b^2\!+\!2 b^3\right) \left(7 a^4\!+\!5 a^3 b\!+\!9 a^2 b^2\!+\!4 a b^3\!+\!2 b^4\right)\\ 
& \hskip -1.2cm t_{0,3}^{(1)}=4 (a^2\!-\!b^2)^2\left(2 a^2\!+\!b^2\right) \left(140 a^6\!+\!41 a^5 b\!+\!193 a^4 b^2\!+\!36 a^3 b^3\!+\!68 a^2 b^4\!+\!4 a b^5\!+\!4 b^6\right)\\ 
& \hskip -1.2cm t_{0,4}^{(1)}=-4 (a\!-\!b)^3 \left(2 a^2\!+\!b^2\right)^2 \left(5 a^2-2 a b\!+\!2 b^2\right) \left(7 a^3-3 a^2 b\!+\!6 a b^2-2 b^3\right)\\ 
\end{split}
\end{equation*}
\begin{equation*}
\begin{split}
& \hskip -1.2cm t_{1,0}^{(1)}= -4 (a^2\!-\!b^2)^2 (68 a^8\!+\!316 a^7 b\!+\!812 a^6 b^2\!+\!1428 a^5 b^3\!+\!1881 a^4 b^4\!+\!1839 a^3 b^5\!+\!1414 a^2 b^6\\& \!+\!710 a b^7\!+\!280 b^8)\\ 
& \hskip -1.2cm t_{1,1}^{(1)}= 16 (a\!+\!b)^2 \left(a^4\!+\!7 a^2 b^2\!+\!b^4\right) \left(34 a^6\!+\!86 a^5 b\!+\!155 a^4 b^2\!+\!179 a^3 b^3\!+\!155 a^2 b^4\!+\!86 a b^5\!+\!34 b^6\right)\\ 
& \hskip -1.2cm t_{1,2}^{(1)}=-24 b(a^2\!-\!b^2)^2\left(2 a^2\!+\!b^2\right) \left(34 a^5\!+\!48 a^4 b\!+\!71 a^3 b^2\!+\!52 a^2 b^3\!+\!30 a b^4\!+\!8 b^5\right)\\ 
& \hskip -1.2cm t_{1,3}^{(1)}=-16 (a\!-\!b)^2 \left(2 a^2\!+\!b^2\right) \left(a^4\!+\!7 a^2 b^2\!+\!b^4\right) \left(17 a^4\!+\!5 a^3 b\!+\!15 a^2 b^2\!+\!2 a b^3\!+\!2 b^4\right) \\ 
& \hskip -1.2cm t_{1,4}^{(1)}= 4 (a^2\!-\!b^2)^2\left(2 a^2\!+\!b^2\right)^2 \left(17 a^4-31 a^3 b\!+\!30 a^2 b^2-22 a b^3\!+\!8 b^4\right)\\ 
\end{split}
\end{equation*}
\begin{equation*}
\begin{split}
& \hskip -1.2cm t_{2,0}^{(1)}= 24 a (a^2\!-\!b^2)\left(a^2\!+\!2 b^2\right) \left(2 a^3\!+\!8 a^2 b\!+\!7 a b^2\!+\!10 b^3\right) \left(2 a^4\!+\!4 a^3 b\!+\!9 a^2 b^2\!+\!5 a b^3\!+\!7 b^4\right)\\ 
& \hskip -1.2cm t_{2,1}^{(1)}= -24 a (a^2\!-\!b^2)^2  \left(a^2\!+\!2 b^2\right) \left(8 a^5\!+\!30 a^4 b\!+\!52 a^3 b^2\!+\!71 a^2 b^3\!+\!48 a b^4\!+\!34 b^5\right)\\ 
& \hskip -1.2cm t_{2,2}^{(1)}= 0\\ 
& \hskip -1.2cm t_{2,3}^{(1)}=24 a (a^2\!-\!b^2)^2 \left(2 a^2\!+\!b^2\right) \left(a^2\!+\!2 b^2\right) \left(4 a^3-3 a^2 b-2 b^3\right)\\ 
& \hskip -1.2cm t_{2,4}^{(1)}=-24 a (a-2 b) (a^2\!-\!b^2)\left(2 a^2\!+\!b^2\right)^2 \left(a^2\!+\!2 b^2\right)\left(a^2-a b\!+\!b^2\right)\\ 
\end{split}
\end{equation*}
\begin{equation*}
\begin{split}
& \hskip -1.2cm t_{3,0}^{(1)}=4 (a^2\!-\!b^2)^2 \left(a^2\!+\!2 b^2\right) \left(4 a^6\!+\!4 a^5 b\!+\!68 a^4 b^2\!+\!36 a^3 b^3\!+\!193 a^2 b^4\!+\!41 a b^5\!+\!140 b^6\right)\\ 
& \hskip -1.2cm t_{3,1}^{(1)}= -16 (a\!-\!b)^2 \left(a^2\!+\!2 b^2\right) \left(a^4\!+\!7 a^2 b^2\!+\!b^4\right) \left(2 a^4\!+\!2 a^3 b\!+\!15 a^2 b^2\!+\!5 a b^3\!+\!17 b^4\right)\\ 
& \hskip -1.2cm t_{3,2}^{(1)}=-24 b (a^2\!-\!b^2)^2 \left(2 a^2\!+\!b^2\right) \left(a^2\!+\!2 b^2\right) \left(2 a^3\!+\!3 a b^2-4 b^3\right) \\ 
& \hskip -1.2cm t_{3,3}^{(1)}= 16 (a\!+\!b)^2 \left(2 a^2\!+\!b^2\right) \left(a^2\!+\!2 b^2\right) \left(a^2-a b\!+\!b^2\right) \left(a^4\!+\!7 a^2 b^2\!+\!b^4\right)\\ 
& \hskip -1.2cm t_{3,4}^{(1)}= -4 (a^2\!-\!b^2)^2  \left(2 a^2\!+\!b^2\right)^2 \left(a^2\!+\!2 b^2\right) \left(a^2-a b\!+\!4 b^2\right)\\ 
\end{split}
\end{equation*}
\begin{equation*}
\begin{split}
& \hskip -1.2cm t_{4,0}^{(1)}=-4 (a\!-\!b)^3 \left(a^2\!+\!2 b^2\right)^2 \left(2 a^2-2 a b\!+\!5 b^2\right) \left(2 a^3-6 a^2 b\!+\!3 a b^2-7 b^3\right)\\ 
& \hskip -1.2cm t_{4,1}^{(1)}= 4 (a^2\!-\!b^2)^2 \left(a^2\!+\!2 b^2\right)^2 \left(8 a^4\!-\!22 a^3 b\!+\!30 a^2 b^2\!-\!31 a b^3\!+\!17 b^4\right)\\ 
& \hskip -1.2cm t_{4,2}^{(1)}=-24 b  (2 a\!-\!b)  (a^2\!-\!b^2) \left(2 a^2\!+\!b^2\right) \left(a^2\!+\!2 b^2\right)^2 \left(a^2-a b\!+\!b^2\right)\\ 
& \hskip -1.2cm t_{4,3}^{(1)}= -4 (a^2\!-\!b^2)^2  \left(2 a^2\!+\!b^2\right) \left(a^2\!+\!2 b^2\right)^2\left(4 a^2\!-\!a b\!+\!b^2\right)\\ 
& \hskip -1.2cm t_{4,4}^{(1)}=4 (a\!+\!b)^4 \left(2 a^2\!+\!b^2\right)^2 \left(a^2\!+\!2 b^2\right)^2 \ ,\\ 
\end{split}
\end{equation*}
together with 
\begin{equation*}
\begin{split}
& \hskip -1.2cm u_{0,0}^{(0)}=-(a\!-\!b)^2 (a\!+\!b) \left(2 a^2\!+\!b^2\right) \left(a^2\!+\!2 b^2\right)\\ 
& \hskip -1.2cm u_{0,1}^{(0)}=4 (a\!-\!b)^2 (a\!+\!b) \left(a^2\!+\!2 b^2\right)^2\\ 
& \hskip -1.2cm u_{0,2}^{(0)}=-(a\!-\!b) \left(a^2\!+\!2 b^2\right) \left(2 a^4\!+\!17 a^2 b^2\!+\!17 b^4\right)\\ 
& \hskip -1.2cm u_{1,0}^{(0)}=4 (a\!-\!b)^2 (a\!+\!b) \left(2 a^2\!+\!b^2\right)^2\\ 
& \hskip -1.2cm u_{1,1}^{(0)}=-8 (a\!+\!b) \left(4 a^2\!+\!a b\!+\!4 b^2\right) \left(a^4\!+\!7 a^2 b^2\!+\!b^4\right)\\ 
& \hskip -1.2cm u_{1,2}^{(0)}=4 (a^2\!-\!b^2) \left(4 a^5\!+\!14 a^4 b\!+\!22 a^3 b^2\!+\!32 a^2 b^3\!+\!19 a b^4\!+\!17 b^5\right)\\ 
& \hskip -1.2cm u_{2,0}^{(0)}= (a\!-\!b) \left(2 a^2\!+\!b^2\right) \left(17 a^4\!+\!17 a^2 b^2\!+\!2 b^4\right)\\ 
& \hskip -1.2cm u_{2,1}^{(0)}=-4 (a^2\!-\!b^2) \left(17 a^5\!+\!19 a^4 b\!+\!32 a^3 b^2\!+\!22 a^2 b^3\!+\!14 a b^4\!+\!4 b^5\right)\\ 
& \hskip -1.2cm u_{2,2}^{(0)}=(a\!+\!b) \left(34 a^6\!+\!76 a^5 b\!+\!137 a^4 b^2\!+\!154 a^3 b^3\!+\!137 a^2 b^4\!+\!76 a b^5\!+\!34 b^6\right)\\ 
\end{split}
\end{equation*}
and 
\begin{equation*}
\begin{split}
& \hskip -1.2cm u_{0,0}^{(1)}=0\\ 
& \hskip -1.2cm u_{0,1}^{(1)}=-12 b (a\!-\!b)^2 (a\!+\!b) \left(a^2\!+\!2 b^2\right)\\ 
& \hskip -1.2cm u_{0,2}^{(1)}=12 b (a\!-\!b)  \left(a^2\!+\!2 b^2\right)^2 \\ 
& \hskip -1.2cm u_{1,0}^{(1)}=-12 a (a\!-\!b)^2 (a\!+\!b) \left(2 a^2\!+\!b^2\right)\\ 
& \hskip -1.2cm u_{1,1}^{(1)}=48 \left(a^2\!+\!a b\!+\!b^2\right) \left(a^4\!+\!7 a^2 b^2\!+\!b^4\right)\\
& \hskip -1.2cm u_{1,2}^{(1)}=-12 (a^2\!-\!b^2)  \left(2 a^4\!+\!6 a^3 b\!+\!11 a^2 b^2\!+\!9 a b^3\!+\!8 b^4\right)\\ 
& \hskip -1.2cm u_{2,0}^{(1)}= -12 a (a\!-\!b) \left(2 a^2\!+\!b^2\right)^2\\ 
& \hskip -1.2cm u_{2,1}^{(1)}=12 (a^2\!-\!b^2) \left(8 a^4\!+\!9 a^3 b\!+\!11 a^2 b^2\!+\!6 a b^3\!+\!2 b^4\right)\\
& \hskip -1.2cm u_{2,2}^{(1)}= -12 (a\!+\!b)^2 \left(a^2\!+\!a b\!+\!b^2\right) \left(4 a^2\!+\!a b\!+\!4 b^2\right)\ .\\ 
\end{split}
\end{equation*}

\section{Calculation of the law for $\phi=n_2/N$ in the local limit} 
Our starting point is the expansion of $X_{s,t}(G(a,\epsilon),G(b,\epsilon))$ when $\epsilon\to 0$, with $G(c,\epsilon)$ as in \eqref{eq:ghscal}.
Expanding the equation \eqref{eq:eqforXst} to increasing orders in $\epsilon$, we deduce the expansion
\begin{equation*}
X_{s,t}(G(a,\epsilon),G(b,\epsilon))=\frac{3 (s+1) (t+1) (s+t+3)}{(s+3) (t+3) (s+t+1)}+\sum_{i\ge 1}\mathfrak{X}_i(s,t,a,b) \epsilon^{2i}\ ,
\end{equation*} 
where the first term is easily deduced from $\epsilon$ the exact expression \eqref{eq:exactXgg} with $g=h=1/12$
and where the $\mathfrak{X}_i$'s are obtained recursively, order by order in $\epsilon^2$.

Expanding \eqref{eq:eqforXst} at order $\epsilon^2$ shows that $\mathfrak{X}_1(s,t,a,b)=0$ and, 
at order $\epsilon^4$, that $\mathfrak{X}_2$ may be written (by linearity) as 
\begin{equation*}
\mathfrak{X}_2(s,t,a,b)=a^4 \xi(s,t)+ b^4 \xi(t,s)\ ,
\end{equation*}
where $\xi(s,t)$ is the solution of some appropriate partial differential equation.
We thus have the following form for the first terms in the expansion:
\begin{equation*}
\begin{split}
X_{s,t}(G(a,\epsilon),G(b,\epsilon))=\frac{3 (s+1) (t+1) (s+t+3)}{(s+3) (t+3) (s+t+1)}& +(a^4 \xi(s,t)+ b^4 \xi(t,s)) \epsilon^4 \\
&+ \mathfrak{X}_3(s,t,a,b) \epsilon^{6}
+O(\epsilon^8)\ ,\\
\end{split}
\end{equation*} 
from which we deduce via \eqref{eq:FsX} the expansion
\begin{equation*}
F(s,G(a,\epsilon),G(b,\epsilon))=\log \left(\frac{s^2 (2 s+3)}{(s+1)^2 (2 s-1)}\right)-(a^4+b^4) \psi(s) \, \epsilon ^4 
+
 \mathfrak{F}_3(s,a,b) \epsilon^{6}
+O(\epsilon^8)\ ,
\end{equation*}
where $\psi(s)$ is directly related to $\xi(s,t)$ and $\mathfrak{F}_3(s,a,b)$ to $\mathfrak{X}_3(s,t,a,b)$. From the very existence of the scaling function,
we may write
\begin{equation}
\begin{split}
\mathcal{F}(S,a,b)&=\lim_{\epsilon\to 0}\,  \frac{1}{\epsilon^3} \, F\left(\left\lfloor S/\epsilon\right\rfloor,G(a,\epsilon),G(b,\epsilon)\right)\\
&=\frac{1}{2\, S^3}+(a^4+b^4)\, \lim_{\epsilon\to 0} \epsilon \, \psi \left(\left\lfloor S/\epsilon\right\rfloor \right)+ 
\lim_{\epsilon\to 0} \left(\epsilon^3  \mathfrak{F}_3\left(\left\lfloor S/\epsilon\right\rfloor,a,b\right)
+O(\epsilon^5)\right)\ . \\
\end{split}
\label{eq:expFSab}
\end{equation}
This is to be compared with the small $S$ expansion of $\mathcal{F}(S,a,b)$, as obtained from its exact expression
of Appendix A, namely
\begin{equation*}
\mathcal{F}(S,a,b)=\frac{1}{2\, S^3}-(a^4+b^4)\, \frac{S}{60}+(a^6+b^6)\, \frac{S^3}{189}+O(S^5)\ .
\end{equation*}
We readily deduce that $\psi(s)\sim s/60$ when $s\to \infty$ and 
\begin{equation*}
 \mathfrak{F}_3(s,a,b) \underset{s \to \infty}{\sim} (a^6+b^6)\, \frac{s^3}{189}
\end{equation*}
(note that the $O(\epsilon^5)$ term in \eqref{eq:expFSab} necessarily leads to an $O(S^5)$ term in $\mathcal{F}(S,a,b)$), which is precisely the announced scaling correspondence.

\begin{figure}
\begin{center}
\includegraphics[width=6cm]{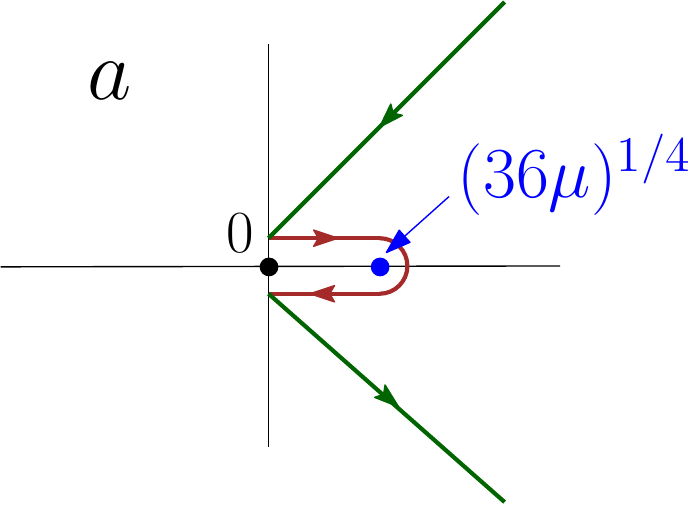}
\end{center}
\caption{Contour of integration in the complex plane for the variable $a$, as inherited from the contour in $g$ around $0$ in
the contour integral \eqref{eq:contour}
via the change of variable $g=G(a,\epsilon)$ when $\epsilon=1/N\to 0$ (see \cite{G17a}).}
\label{fig:contour}
\end{figure}
We may now estimate $[g^N]F(s,g,g\, e^{\mu/N})$.This quantity is obtained by a contour integral around $g=0$, namely
\begin{equation}
\frac{1}{2\rm{i}\pi} \oint\frac{dg}{g^{N+1}} F(s,g,g\, e^{\frac{\mu}{N}})
\label{eq:contour}
\end{equation}
and, at large $N$, we may change variable from $g$ to $a$ by taking $g=G(a,\epsilon)$ with
\begin{equation*}
\hskip 5.cm \epsilon=\frac{1}{N^{1/4}}\ .
\end{equation*}
Setting $h=g\, e^{\frac{\mu}{N}} =g\, (1+ \mu \epsilon^4)$ then amounts to choosing
\begin{equation*}
\frac{-b^4}{36}=\frac{-a^4}{36}+\mu+O\left(\frac{1}{N}\right) \quad \Leftrightarrow \quad b^4=a^4-36\mu +O\left(\frac{1}{N}\right)\ .
\end{equation*}
Using $dg=-(1/12)a^3/(9 N)$ and $g^{N+1}\sim (1/12)^{N+1}\, e^{-a^4/36}$, we eventually arrive at
\begin{equation*}
\begin{split}
[g^N]F(s,g,g\, e^{\mu/N}) & =  \frac{1}{2\rm{i}\pi} \frac{12^N}{N} \int_{\mathcal{C}_\mu} da\,  \frac{-a^3}{9}
\, e^{a^4/36} \Big\{\log \left(\frac{s^2 (2 s+3)}{(s+1)^2 (2 s-1)}\right)\\
& -\frac{1}{N}\, (2a^4-36\mu ) \psi(s) 
+\frac{1}{N^{3/2}}
 \mathfrak{F}_3(s,a,(a^4-36 \mu)^{1/4})
+O\left(\frac{1}{N^2}\right)
\Big\}\\
\end{split}
\end{equation*}
with some appropriate integration contour $\mathcal{C}_\mu$ in the complex plane. As explained in \cite{G17a} and illustrated in Figure~\ref{fig:contour}, this contour
is made of a first part $\mathcal{C}_{\mu,1}$ consisting of two half straight lines at $\pm 45^\circ$ meeting at the origin, and a part $\mathcal{C}_{\mu,2}$ 
which makes a back and forth excursion from $0$ to $(36\mu)^{1/4}$. 
Both the constant (i.e.\ independent of $a$) terms and the $a^4$ term in-between the curly brackets lead to integrals along this contour which \emph{vanish identically} by symmetry\footnote{This vanishing holds for any finite $s$, i.e.\ even before taking the $s\to \infty$ limit.}, 
so that
\begin{equation*}
[g^N]F(s,g,g\, e^{\mu/N}) \underset{N \to \infty}{\sim}  \frac{1}{2\rm{i}\pi} \frac{12^N}{N^{5/2}} \int_{\mathcal{C}_\mu} da\,  \frac{-a^3}{9}
\, e^{a^4/36} \mathfrak{F}_3(s,a,(a^4-36 \mu)^{1/4})
\end{equation*}
with a right hand side which behaves at large $s$ as
\begin{equation*}
\hskip -1.2cm  \frac{1}{2\rm{i}\pi}\int_{\mathcal{C}_\mu} da\,  \frac{-a^3}{9}
\, e^{a^4/36} \mathfrak{F}_3(s,a,(b^4\!-\!36 \mu)^{1/4})\underset{s \to \infty}{\sim}\frac{s^3}{189} \frac{1}{2\rm{i}\pi} \int_{\mathcal{C}_\mu} da\,  \frac{-a^3}{9}
\, e^{a^4/36} (a^6+(a^4\!-\!36 \mu)^{3/2})\ .
\end{equation*}
The contribution of the two terms in this latter integral were computed in \cite{G17a}, namely
\begin{equation*}
\begin{split}
&\frac{1}{2\rm{i}\pi} \int_{\mathcal{C}_\mu} da\,  \frac{-a^3}{9}
\, e^{a^4/36} a^6 = \frac{432}{\pi}\, \int_0^\infty dt\, t^4e^{-t^2}=\frac{162}{\sqrt{\pi}}\ , \\
&
\frac{1}{2\rm{i}\pi} \int_{\mathcal{C}_\mu} da\,  \frac{-a^3}{9}
\, e^{a^4/36} (a^4\!-\!36 \mu)^{3/2}= \frac{432\, e^{\mu}}{\pi}\, \int_0^\infty dt\, t^4e^{-t^2}=\frac{162}{\sqrt{\pi}}\, e^{\mu}\ .\\
\end{split}
\end{equation*}
This yields eventually
\begin{equation*}
[g^N]F(s,g,g\, e^{\mu/N}) \underset{N \to \infty}{\sim}  \frac{12^N}{\sqrt{\pi}\, N^{5/2}} 
\left(\frac{6}{7} (1+e^{\mu})\, s^3 + O(s^2)\right)
\end{equation*}
at large $s$.
At $\mu=0$, we recover the estimate \eqref{eq:larges}. Taking the appropriate ratio, we arrive immediately at the desired result \eqref{eq:trivial}.

\section{Estimate of $\mathfrak{P}(V)$ at small $V$} 
Since the quantity 
\begin{equation*}
\hat{\mathfrak{P}}(\sigma)\equiv E\left[e^{-\sigma\, V}\right]
\end{equation*}
is the Laplace transform of the probability distribution $\mathfrak{P}(V)$, with $V$ taking its values in $[0,\infty)$, 
it has no singularity for real non-negative $\sigma$. Its inverse Laplace transform, the probability distribution $\mathfrak{P}(V)$ itself, may thus be
obtained via
\begin{equation*}
\mathfrak{P}(V)= \frac{1}{2\rm{i}\pi}\int_{\gamma-\rm{i}\infty}^{\gamma+\rm{i}\infty}d\sigma\, e^{\sigma\, V}\, \hat{\mathfrak{P}}(\sigma)
\end{equation*}
for any real non-negative $\gamma$. 
At small $V$, this integral may be evaluated via a saddle point approximation as follows. For large $\sigma$, we have
the asymptotic equivalence
\begin{equation*}
\hat{\mathfrak{P}}(\sigma) \underset{\sigma \to \infty}{\sim} \frac{9}{2}(3\sqrt{2}-4)\, e^{-\sqrt{6}\, \sigma^{1/4}}
 \end{equation*}
and the integral is dominated by its saddle point $\sigma^*$ given by
\begin{equation*}
\frac{d\ }{d\sigma}\left( \sigma\, V -\sqrt{6}\, \sigma^{1/4}\right)\Big\vert_{\sigma=\sigma^*}=0\quad
\Leftrightarrow \quad \sigma^*= \frac{3^{2/3}}{4 V^{4/3}}\ .
\end{equation*}
The use of the asymptotic equivalent above for $\hat{\mathfrak{P}}(\sigma)$ is fully consistent if $\sigma^*$ becomes large,
i.e. when $V$ itself becomes \emph{small}. Setting
\begin{equation*}
\sigma= \sigma^*+\rm{i}\, \eta
\end{equation*}
in the integral, with $\eta$ real (i.e.\ choosing implicitly $\gamma=\sigma^*$), we may use the expansion
\begin{equation*}
 \sigma\, V -\sqrt{6}\, \sigma^{1/4}= -\frac{3^{5/3}}{4\, V^{1/3}}-\frac{3^{1/3}}{2}\,  V^{7/3}\, \eta ^2+O\left(\eta ^3\right)
\end{equation*}
to write
\begin{equation*}
\begin{split}
\mathfrak{P}(V)\underset{V \to 0}{\sim} & \ \frac{9}{2}(3\sqrt{2}-4)\,  e^{-\frac{3^{5/3}}{4\, V^{1/3}}} \frac{1}{2 \pi}\int_{-\infty}^{\infty}d\eta\, e^{-\frac{3^{1/3}}{2}\,  V^{7/3}\, \eta ^2}\\
& = \frac{3^{11/6} (3 -2 \sqrt{2})}{2 \sqrt{\pi }}\, \frac{1}{V^{7/6}}\, e^{-\frac{3^{5/3}}{4}\frac{1}{V^{1/3}}}\ .
\\
\end{split}
\end{equation*}

\bibliographystyle{plain}
\bibliography{voronoiinfinite}
\end{document}